\documentclass{amsart}
\usepackage{fullpage}
\usepackage{amsmath}
\usepackage{amsthm}
\usepackage{amssymb}
\usepackage{tikz}
\usepackage{tikz-cd}
\usepackage{mathdots}
\usepackage{mathtools}
\usepackage{enumitem}
\usepackage{thmtools}
\usepackage[hypertexnames=false]{hyperref}
\usepackage{cleveref}
\usepackage{xcolor}
\hypersetup{
    colorlinks = true,
    linkcolor={red!50!black},
    citecolor={blue!50!black},
    urlcolor={blue!80!black}
}
\usepackage[T1]{fontenc}
\usepackage[utf8]{inputenc}

\declaretheorem[name=Theorem, refname={Theorem,Theorems}]{maintheorem}

\numberwithin{equation}{section}
\declaretheorem[name=Theorem, sibling=equation, refname={Theorem,Theorems}]{theorem}
\declaretheorem[name=Proposition, sibling=equation]{proposition}
\declaretheorem[name=Lemma, sibling=equation, refname={Lemma,Lemmas}]{lemma}
\declaretheorem[name=Corollary, sibling=equation, refname={Corollary,Corollaries}]{corollary}
\declaretheorem[name=Conjecture, sibling=equation]{conjecture}
\declaretheorem[name=Guess, sibling=equation]{guess}
\declaretheorem[name=Sublemma, numbered=no]{sublemma}
\theoremstyle{definition}
\newtheorem{remark}[equation]{Remark}
\newtheorem{definition}[equation]{Definition}
\newtheorem{question}[equation]{Question}
\setlist[enumerate]{label=\textup{(}\alph*\textup{)}}
\newcommand{\ring}[1]{\mathbb{#1}}
\newcommand{\Z}{\ring{Z}}
\newcommand{\Fp}{\ring{F}\!_{p}}
\newcommand{\ftwo}{\ring{F}\!_{2}}
\DeclareMathOperator{\Hom}{Hom}
\DeclareMathOperator{\Ext}{Ext}
\newcommand{\tildePP}[1]{\widetilde{\mathcal{P}}^{#1}}
\newcommand{\PP}[1]{\mathcal{P}^{#1}}
\newcommand{\pst}[2]{P^{#1}_{#2}}
\newcommand{\modmod}{/ \! /}
\newcommand{\massey}[1]{\langle #1 \rangle}
\newcommand{\cat}[1]{\ensuremath{\mathsf{#1}}}
\newcommand{\stable}[1]{\cat{Stable}(\dual{#1})}
\newcommand{\Bdd}{\cat{Bdd}}
\newcommand{\nBdd}{\cat{nBdd}}
\newcommand{\cotensor}[1]{\,\Box_{#1}}
\newcommand{\exterior}{\Lambda}
\newcommand{\restrict}[1]{\!\downarrow_{#1}}
\newcommand{\Smash}{\wedge}
\newcommand{\alg}[1]{#1}
\newcommand{\plaindual}[1]{#1^{\smash{\scalebox{.7}[1.4]{\rotatebox{90}{\textnormal{\guilsinglleft}}}}}}
\newcommand{\dual}[1]{\plaindual{\alg{#1}}}
\newcommand{\Binom}{\alg{Bi}}
\renewcommand{\P}{\alg{P}}
\newcommand{\D}{\alg{D}}
\newcommand{\A}{\alg{A}}
\renewcommand{\H}[1]{H\!\alg{#1}}
\DeclarePairedDelimiter\floor{\lfloor}{\rfloor}
\DeclarePairedDelimiter\bousfield{\langle}{\rangle}
\begin{document}

\title{Toward Quillen stratification for \\
the odd primary Steenrod algebra}
\author{John H. Palmieri}

\begin{abstract}
Fix an odd prime $p$ and let $\A$ be the mod $p$ Steenrod
algebra. There is a normal sub-Hopf algebra $\D$ of $\A$, and we show
that the cohomology of $\A$ is $F$-isomorphic to the subalgebra of the
cohomology of $\D$ which is invariant under the $\A$-action. We also
show that the sub-Hopf algebras of $\A$ with bounded profile function
detect nilpotence in cohomology.
\end{abstract}

\subjclass{Primary 55S10; 
Secondary 55T15, 
16T05, 
18G15} 

\maketitle


\section{Introduction}

Throughout, $p$ will be an odd prime number.

Let $\A$ be the mod $p$ Steenrod algebra. The cohomology of $\A$,
$\Ext_{\A}^{*}(\Fp, \Fp)$, is the $E_{2}$-page of the classical Adams
spectral sequence converging to the $p$-component of the stable
homotopy groups of spheres, an object of great interest in stable
homotopy theory. The goal of this paper is to try to understand the
cohomology of $\A$ along the lines of the so-called Quillen
stratification for the mod 2 Steenrod algebra, as presented in
\cite{palmieri;quillen} and \cite[Chapter
5]{palmieri;steenrod-stable}.

As Milnor \cite{milnor;steenrod} tells us, the Steenrod algebra is the
Hopf algebra with graded dual
\[
\dual{A} = \Fp [\xi_{1}, \xi_{2}, \xi_{3}, \dots] \otimes 
\exterior (\tau_{0}, \tau_{1}, \tau_{2}, \dots),
\]
graded by $\deg \xi_{n} = 2(p^{n}-1)$ and $\deg \tau_{n} =
2p^{n}-1$, with coproduct
\begin{align*}
\xi_{n} &\mapsto \sum_{i=0}^{n} \xi_{n-i}^{p^{i}} \otimes \xi_{i}, \\
\tau_{n} &\mapsto \tau_{n} \otimes 1 + 
\sum_{i=0}^{n} \xi_{n-i}^{p^{i}} \otimes \tau_{i},
\end{align*}
where $\xi_{0}=1$ when it appears in these formulas.

To understand a Hopf algebra like $\A$, a good first step is to
understand its sub-Hopf algebras. Any sub-Hopf algebra $\alg{C}$ of
$\A$ is dual to a quotient Hopf algebra $\dual{C}$ of $\dual{A}$, and
by a theorem of Adams and Margolis \cite{adams-margolis;alg}, those
quotients must have the form
\[
\dual{C} = \dual{A} / (\xi_{1}^{p^{c_{1}}},
\xi_{2}^{p^{c_{2}}}, \xi_{3}^{p^{c_{3}}}, \dots; \tau_{0}^{e_{0}},
\tau_{1}^{e_{1}}, \tau_{2}^{e_{2}}, \dots)
\]
for some exponents $((c_{1}, c_{2}, c_{3}, \dots); (e_{0},
e_{1}, e_{2}, \dots))$, where $c_{i}$ is either a non-negative integer
or $\infty$ and $e_{i} \in \{1,2 \}$. The list of exponents is
called the \emph{profile function} for the Hopf algebra.

Adams and Margolis characterized precisely which such lists correspond
to sub-Hopf algebras, thus providing a complete classification. We
review the details in the next section and for now just mention one
case: we let $\D$ be the sub-Hopf algebra of $\A$ with profile
function
\[
((1, 2, 3, \dots, n, \dots); (2, 2, 2, \dots).
\]
This algebra is normal in $\A$, and more generally, for any sub-Hopf
algebra $\alg{C}$ of $\A$, $\alg{C} \cap \D$ is normal in $\alg{C}$.
The inclusion map $i: \alg{C} \cap \D \to \alg{C}$ induces a
restriction map $i^{*}: \Ext_{\alg{C}}^{*}(\Fp, \Fp) \to \Ext_{\alg{C}
\cap \D}^{*}(\Fp, \Fp)$. Because $\alg{C} \cap \D$ is normal in
$\alg{C}$, the restriction map factors through the subalgebra of
invariants under the action of the quotient Hopf algebra $\alg{C}
\modmod (\alg{C} \cap \D) := \alg{C} \otimes_{\alg{C} \cap \D} \Fp$:
\[
\begin{tikzcd}
\Ext_{\alg{C}}^{*}(\Fp, \Fp) \arrow[r, "i^{*}"] \arrow[rd,
"\overline{\imath}^{*}"] & \Ext_{\alg{C} \cap \D}^{*}(\Fp, \Fp) \\
& \left( \Ext_{\alg{C} \cap \D}^{*}(\Fp,
\Fp)\right)^{\alg{C} \modmod \alg{C} \cap \D} \arrow[u, hook]
\end{tikzcd}
\]
The map $\overline{\imath}^{*}$ will make an appearance in
\autoref{thm-f-iso} below.

Quillen stratification involves understanding cohomology algebras up
to $F$-isomorphism, so recall:

\begin{definition}
A map $\varphi$ of graded
commutative algebras over a field of characteristic $p$ is an
\emph{$F$-isomorphism} if
\begin{itemize}
\item every element in the kernel of $\varphi$ is nilpotent, and
\item for every $y$ in the codomain of $\varphi$, there is an $n$ so
that $y^{p^{n}}$ is in the image of $\varphi$.
\end{itemize}
\end{definition}

We can strengthen the first part of the definition of $F$-isomorphism
to echo the nilpotence theorem of Devinatz, Hopkins, and Smith. To do
that, we work in a category which behaves like the classical stable
homotopy category: given a sub-Hopf algebra $\alg{C}$ of $\A$, let
$\stable{C}$ be the category as defined in \cite[Section
9.5]{axiomatic} and examined in detail in
\cite{palmieri;steenrod-stable}. The objects of $\stable{C}$ are
cochain complexes of injective graded left comodules over $\dual{C}$,
and the morphisms are cochain homotopy classes of maps. This category
is a stable homotopy category in the sense of \cite{axiomatic}, and we
may refer to its objects as ``spectra.''  In this category, the sphere
object $S$ is an injective resolution of $\Fp$, and if $X$ is any
object, then we may write $\pi_{*}(X)$ for the chain homotopy classes
of maps (of any degree) from $S$ to $X$. If $X$ is the injective
resolution of a comodule $M$, then $\pi_{*}(X) \cong
\Ext_{\alg{C}}^{*}(\Fp, M)$. (Note that if $M$ is a
$\dual{C}$-comodule, then $M$ can also be viewed as a
$\alg{C}$-module, and so Ext is computed in the category of
$\alg{C}$-modules.) Also in this category, any sub-Hopf algebra
$\alg{B}$ of $\alg{C}$ defines a ring spectrum called $\H{B}$: this is
an injective resolution of $\dual{C} \cotensor{\dual{B}} \Fp \cong
\plaindual{(\alg{C} \otimes_{\alg{B}} \Fp)}$, and its coefficient ring
is $\pi_{*}\H{B} = \H{B}_{*} \cong \Ext_{\alg{B}}^{*}(\Fp, \Fp)$. If
$X$ is the injective resolution of a $\dual{C}$-comodule $M$, then by
restriction $M$ has the structure of a $B$-module, and we have
$\H{B}_{*}(X) \cong \Ext_{\alg{B}}^{*}(\Fp, M)$.

\begin{definition}
Given a Hopf algebra $\alg{C}$,
we say that a collection $\cat{S}$ of sub-Hopf algebras of $\alg{C}$
\emph{detects nilpotence} in $\stable{C}$ if the following
conditions hold:
\begin{enumerate}
\item Fix a ring spectrum $R$ and $\alpha \in \pi_{*}(R)$. Then
$\alpha$ is nilpotent if $\H{B}_{*}\alpha = 0$ for all $\alg{B} \in
\cat{S}$.
\item Fix a finite spectrum $Y$ and a self-map $f: \Sigma^{m} Y \to
Y$. Then $f$ is nilpotent if and only if $\H{B}_{*}f$ is nilpotent for
all $\alg{B} \in \cat{S}$.
\item Fix a finite spectrum $F$, an arbitrary spectrum $X$, and a map
$f: F \to X$. Then $f$ is smash-nilpotent if $1_{\H{B}} \Smash f: \H{B}
\Smash F \to \H{B} \Smash X$ is zero for all $\alg{B} \in \cat{S}$.
\end{enumerate}
\end{definition}

We have three main results in this paper; here are the first
two.

\begin{maintheorem}\label{thm-nilpotence}
Fix an odd prime $p$ and suppose that $\alg{C}$ is a sub-Hopf algebra
of $\A$.  The sub-Hopf algebra $\alg{C} \cap \D$ detects nilpotence in
$\stable{C}$.
\end{maintheorem}

\begin{maintheorem}\label{thm-f-iso}
Fix an odd prime $p$ and suppose that $\alg{C}$ is a sub-Hopf algebra of
$\A$. The map
\[
\overline{\imath}^{*}: \Ext_{\alg{C}}^{*}(\Fp, \Fp) \to \left( \Ext_{\alg{C}
\cap \D}^{*}(\Fp, \Fp)\right)^{\alg{C} \modmod \alg{C} \cap \D}
\]
is an $F$-isomorphism. 
\end{maintheorem}

For example if $\alg{C}=\A$, the first theorem tells us that every element in
the kernel of the restriction map $i^{*}: \Ext_{\A}^{*}(\Fp, \Fp) \to
\Ext_{\D}^{*}(\Fp, \Fp)$ is nilpotent, and the second theorem tells us
that the map
\[
\overline{\imath}^{*}: \Ext_{\A}^{*}(\Fp, \Fp) \to \left(
\Ext_{\D}^{*}(\Fp, \Fp)\right)^{\A \modmod \D}
\]
is an $F$-isomorphism.

These theorems also hold at the prime 2 \cite{palmieri;quillen}, and
in that setting one can use them to go further, in particular to
obtain results involving ``elementary'' sub-Hopf algebras. A
cocommutative Hopf algebra $\alg{E}$ over a field $k$ of
characteristic $p>0$ is called \emph{elementary} if it is commutative
and satisfies $(I\!\alg{E})^{p}=0$, where $I\!\alg{E}$ is the
augmentation ideal of $\alg{E}$. At the prime 2, this means that
$\alg{E}$ is isomorphic as an algebra to $k[x_{1}, x_{2},
\dots]/(x_{i}^{2})$, whereas at odd primes $\alg{E}$ is isomorphic to
a tensor product of a polynomial algebra truncated at height $p$ and
an exterior algebra:
\[
\alg{E} \cong_{\cat{alg}} k[x_{1}, x_{2}, \dots]/(x_{i}^{p}) \otimes
\exterior (y_{1}, y_{2}, \dots).
\]
Up to $F$-isomorphism the cohomology of an elementary Hopf algebra is
a polynomial algebra. For the mod 2 Steenrod algebra $\A$, we have
this:

\begin{theorem}[Quillen stratification for the mod $2$ Steenrod
algebra \cite{palmieri;quillen}]\label{thm-2-elementary}
We work at the prime 2. Fix a sub-Hopf algebra $\alg{C}$ of the mod
2 Steenrod algebra.
\begin{enumerate}
\item The elementary sub-Hopf algebras of $\alg{C}$ detect nilpotence.
\item The map 
\[
\Ext_{\alg{C}}^{*}(\ftwo, \ftwo) \to \left( \varprojlim_{\alg{E}}
\Ext_{\alg{E}}^{*}(\ftwo, \ftwo) \right)^{\alg{C}}
\]
is an $F$-isomorphism, where the limit is taken over the category
whose objects are the elementary sub-Hopf algebras of $\alg{C}$ and
whose morphisms are inclusions.
\end{enumerate}
\end{theorem}

This theorem deserves to be called Quillen stratification for the mod
2 Steenrod algebra. If $k$ is a field of characteristic $p$, Quillen
stratification for a finite group $G$ provides a description of
$H^{*}(G;k)$: there is an $F$-isomorphism
\[
H^{*}(G;k) \to \varprojlim_{\cat{Elem}(G)} H^{*}(E;k),
\]
where $\cat{Elem}(G)$ is the category whose objects $E$ are elementary
abelian $p$-subgroups of $G$ and whose morphisms are generated by
inclusions and conjugations by elements of $G$. For elementary sub-Hopf
algebras of the Steenrod algebra, it turns out that one can separate
the inclusions and the conjugations, so the inclusions combine to give
$\varprojlim_{\alg{E}} \Ext_{\alg{E}}^{*}(\ftwo, \ftwo)$, while the
analog of taking the limit under conjugations is to take the
invariants.

Wilkerson showed \cite{wilkerson} that the elementary sub-Hopf
algebras of the odd primary Steenrod algebra do not detect nilpotence,
so the precise analog of this at odd primes cannot hold. We can,
though, provide a partial replacement. We describe the elementary
sub-Hopf algebras of the odd primary Steenrod algebra in the next
section; the relevant parts here are that (a) each elementary $\alg{B}
\subseteq \A$ has a \emph{bounded profile function}: the profile
function has the form $((b_{1}, b_{2}, \dots); (e_{0}, e_{1}, \dots))$
such that there is a number $M$ with each $b_{i} \leq M$; and (b) each
elementary $\alg{B} \subseteq \A$ is contained in $\D$; indeed, $\D$
is the smallest Hopf algebra containing every elementary sub-Hopf
algebra. There are, though, sub-Hopf algebras of $\D$ with bounded
profile functions which are not elementary, say for example the one
given by the profile function $((1, 2, 3, 3, 3, 3, \dots); (2, 2, 2,
\dots))$. The following partial replacement for
\autoref{thm-2-elementary} is the third main result of this paper.

\begin{maintheorem}\label{thm-bdd-nilpotence}
Fix a sub-Hopf algebra $\alg{C} \subseteq \A$. Let $\Bdd(\alg{C} \cap
\D)$ be the collection of sub-Hopf algebras of $\alg{C} \cap \D$ with
bounded profile function. Then $\Bdd(\alg{C} \cap \D)$ detects
nilpotence in $\stable{C}$.
\end{maintheorem}

One might also conjecture that there is an $F$-isomorphism
\[
H^{*}(\alg{C}) \to \left( \varprojlim_{\Bdd (\alg{B} \cap \D)}
H^{*}(\alg{B}) \right)^{\alg{C}},
\]
where $\Bdd (\alg{B} \cap \D)$ is made into a category with the
morphisms being inclusions. This conjecture is discussed in
\autoref{sec-bounded}. Also see \autoref{conj-f-iso-quasi}
for a stronger statement, a more direct analog of the $p=2$ theorem.

The structure of the paper: in \autoref{sec-prelim} we provide
background, notation, and a few key lemmas. In
\autoref{sec-reduce} we prove a theorem that provides a framework
within the Steenrod algebra for proving results like
\cref{thm-nilpotence,thm-f-iso}. We then use that
theorem in \autoref{sec-binom} to prove an intermediate result: a
certain sub-Hopf algebra $\Binom$ (where $\D \subseteq \Binom
\subseteq \A$) detects nilpotence and gives an $F$-isomorphism, and we
use that in \autoref{sec-D} to show that $\D$ detects nilpotence
and gives an $F$-isomorphism, yielding our first two main theorems. In
\autoref{sec-bounded} we prove \autoref{thm-bdd-nilpotence},
and in \autoref{sec-discussion} we discuss some further
directions. A key computation used in \autoref{sec-binom} comes
from unpublished work of Burklund, Levy, and Pstr\k{a}gowski; it is
valid at all primes, but in an appendix we give an independent proof
for the $p=3$ case.


\section{Preliminaries}
\label{sec-prelim}

As we have said, $\A$ denotes the mod $p$ Steenrod algebra with
$p>2$. Throughout, (co)modules will be understood to be left
(co)modules.


\subsection{Hopf algebras related to \texorpdfstring{$\A$}{A}}\label{subsec-hopf}

First we recall the Adams-Margolis theorem classifying sub-Hopf
algebras of $\A$ \cite{adams-margolis;alg}. We consider a
sub-coalgebra $\alg{B}$ of $\A$ whose graded dual is the quotient
algebra of $\dual{A}$ of the form
\[
\dual{B} = \dual{A} / (\xi_{1}^{p^{b_{1}}},
\xi_{2}^{p^{b_{2}}}, \xi_{3}^{p^{b_{3}}}, \dots; \tau_{0}^{e_{0}},
\tau_{1}^{e_{1}}, \tau_{2}^{e_{2}}, \dots)
\]
for some exponents $((b_{1}, b_{2}, b_{3}, \dots); (e_{0}, e_{1},
e_{2}, \dots))$, where $b_{i} \in \{0, 1, 2, \dots, \infty \}$ and
$e_{i} \in \{1,2 \}$. (If some $b_{i} = \infty$, then do not mod out
by any power of $\xi_{i}$.) We have the following conditions:
\begin{itemize}
\item for all $i, j \geq 1$, either $b_{i} \leq j + b_{i+j}$ or $b_{j}
\leq b_{i+j}$, and
\item for all $i \geq 1$ and $j \geq 0$, if $e_{i+j}=1$, then either
$b_{i} \leq j$ or $e_{j} = 1$.
\end{itemize}

\begin{theorem}[Theorem 2.4 in \cite{adams-margolis;alg}]
\label{thm-adams-margolis}
If the sub-coalgebra $\alg{B} \subseteq \A$ satisfies these
conditions, then it is a sub-Hopf algebra, and conversely, every
sub-Hopf algebra of $\A$ is a sub-coalgebra of this form satisfying
these conditions.
\end{theorem}

We also need to understand when a sub-Hopf algebra of $\A$ is
normal. If $\alg{H}$ is a cocommutative Hopf algebra over $\Fp$, then
a sub-Hopf algebra $\alg{B} \subseteq \alg{H}$ is \emph{normal} if the
left and right ideals generated by $\alg{B}$ are equal: $\alg{H} \cdot
I\! \alg{B} = I\! \alg{B} \cdot \alg{H}$, where $I\! \alg{B} = \ker (\alg{B}
\xrightarrow{\varepsilon} \Fp)$ is the augmentation ideal of $\alg{B}$. It is well
documented, for example in \cite{margolis}, which sub-Hopf algebras of
the mod 2 Steenrod algebra are normal: those with non-decreasing
profile function. The author has not been able to find a good
reference for the odd primary case, though.

\begin{proposition}
Suppose that $\alg{B} \subseteq \A$ is a sub-Hopf algebra with
profile function
\[
((b_{1}, b_{2}, b_{3}, \dots); (e_{0}, e_{1}, e_{2}, \dots)).
\]
Then $\alg{B}$ is normal in $\A$ if and only if
\begin{enumerate}
\item the sequence $(b_{i})$ is non-decreasing, 
\item the sequence $(e_{i})$ is non-decreasing, and
\item if $e_{n}=1$ for some $n$, then $b_{i}=0$ for $i \leq n$.
\end{enumerate}
\end{proposition}

\begin{proof}
The sub-Hopf algebra $\alg{B} \subseteq \A$ is normal if and only if
there is an extension of Hopf algebras
\[
\alg{B} \to \A \to \A \modmod \alg{B},
\]
and this happens if and only if there is an extension of Hopf algebras 
\[
\dual{A} \cotensor{\dual{B}} \Fp \to \dual{A} \to \dual{B}.
\]
We have
\[
\dual{B} = \dual{A} / (\xi_{1}^{p^{b_{1}}},
\xi_{2}^{p^{b_{2}}}, \xi_{3}^{p^{b_{3}}}, \dots; \tau_{0}^{e_{0}},
\tau_{1}^{e_{1}}, \tau_{2}^{e_{2}}, \dots),
\]
and so
\[
\dual{A} \cotensor{\dual{B}} \Fp = \Fp [\xi_{1}^{p^{b_{1}}}, \xi_{2}^{p^{b_{2}}},
\xi_{3}^{p^{b_{3}}}, \dots] \otimes \exterior (\tau_{n} : e_{n} = 1).
\]
Let $\dual{C} = \dual{A} \cotensor{\dual{B}} \Fp$.  It is clear that
if the conditions hold, then $\dual{C}$ is a sub-Hopf algebra: if the
conditions hold, then $\alg{B}$ is normal. We need to check the
converse, so assume that $\dual{C}$ is a sub-Hopf algebra. Since the
coproduct on $\xi_{n}^{p^{b_{n}}}$ contains terms
$\xi_{n-i}^{p^{b_{n}+i}} \otimes \xi_{i}^{p^{b_{n}}}$, we conclude
that $b_{i} \leq b_{n}$ if $i < n$. Similarly, the coproduct on
$\tau_{n}$ contains terms $\xi_{n-i}^{p^{i}} \otimes \tau_{i}$, so if
$e_{n}=1$ (meaning that $\tau_{n}$ is a generator of $\dual{C}$), then
$e_{i}=1$ when $i < n$, and $\xi_{n}$ must lie in $\dual{C}$: the
exponent $b_{n}$ must be 0. Since $b_{i} \leq b_{n}$ for $i \leq n$,
we conclude that $b_{i} = 0$ for all $i \leq n$.
\end{proof}

In practice we will typically cite the conditions to conclude that a
sub-Hopf algebra is normal, and we will also typically be in the
situation where each $e_{i}=2$; or somewhat more generally we will
have Hopf algebras $\alg{B}' \subseteq \alg{B}$ which agree on the
exterior parts of their profile functions, and we will want to know
when $\alg{B}'$ is normal in $\alg{B}$. 

We define some sub- and quotient Hopf algebras of $\A$.
\begin{itemize}
\item (This one is only used to help define some of the others in this
list.) Write $\dual{\exterior}$ for the exterior quotient Hopf algebra
of $\dual{A}$: $\dual{\exterior} := \exterior (\tau_{0}, \tau_{1},
\tau_{2}, \dots)$.  Its graded dual is $\alg{\exterior} := \exterior
(Q_{0}, Q_{1}, Q_{2}, \dots) \subseteq \A$, which is normal in $\A$.
\item We have already introduced $\D \subseteq \A$, defined by
$\dual{D} \cong \dual{A} / (\xi_{1}^{p}, \xi_{2}^{p^{2}}, \dots,
\xi_{n}^{p^{n}}, \dots)$.  Its profile function is therefore
\[
((1,2,3, \dots, n, \dots); (2, 2, 2, \dots)).
\]
This is called ``$\!\D$'' because the polynomial part of the profile
function traces the ``diagonal'' line $n \mapsto n$. $\D$ is normal in
$\A$. As noted above, $\D$ is the algebra generated by the elementary
sub-Hopf algebras of $\A$. It is also the algebra generated by the
elements $\pst{s}{t}$ with $s<t$ and $Q_{n}$ for all $n$; these
elements define the Margolis homology groups which detect projectivity
for $\A$-modules. So $\D$ is a natural algebra to consider when
investigating cohomological properties of $\A$.
\item Define $\Binom \subseteq \A$ by
$\dual{\Binom} \cong \dual{A} / (\xi_{1}^{p}, \xi_{2}^{p^{2}},
\xi_{3}^{p^{4}}, \dots, \xi_{n}^{p^{\binom{n}{2}+1}}, \dots)$,
so its profile function is
\[
((1, 2, 4, 7, \dots, \binom{n}{2}+1, \dots); (2, 2, 2, \dots)).
\]
In case you're curious, ``$\!\Binom$'' stands for ``binomial.''
$\Binom$ is normal in $\A$. This algebra will be used as a sort of
catalyst: it will appear in some auxiliary results, in particular
\cref{cor-binom,cor-D}, which we use to prove our
main results.
\item Define $\P$ to be dual to the polynomial part of $\dual{A}$,
$\dual{P} = \Fp [\xi_{1}, \xi_{2}, \xi_{3}, \dots]$.
Since $\dual{P}$ is a sub-Hopf algebra of $\dual{A}$, we see that $\P$
is naturally a quotient Hopf algebra of $\A$: $\P \cong \A \modmod
\alg{\exterior}$.  One can also view $\P$ as a sub-Hopf algebra of
$\A$: use the profile function
\[
((\infty, \infty, \infty, \dots); (1, 1, 1, \dots)).
\]
As a consequence, $\A$ can be viewed a semi-direct product of $\P$ and
$\alg{\exterior}$.  Because the composite $\P \to \A \to \P$ of
inclusion followed by projection is the identity, we see that the
cohomology of $\P$ splits off of the cohomology of $\A$. Somewhat like
$\Binom$, $\P$ is used as a catalyst, but it appears pretty naturally
because of the properties just described.
\item Define $\A_{\geq n} \subseteq \A$ by $\dual{A}_{\geq n} \cong
\dual{A} / (\xi_{1}, \dots, \xi_{n-1}, \tau_{0}, \dots, \tau_{n-1})$.
Its profile function is therefore
\[
((\underbrace{0, \dots, 0}_{n-1}, \infty, \infty, \infty, \dots); 
(\underbrace{1, \dots, 1}_{n}, 2, 2, 2, \dots)).
\]
$\A_{\geq n}$ is normal in $\A$.
\item Let $\A_{< n}$ be the quotient Hopf algebra $\A \modmod
\A_{\geq n}$. Thus we have 
\[
\dual{A}_{< n} \cong \Fp [\xi_{1}, \dots, \xi_{n-1}] \otimes 
\exterior (\tau_{0}, \dots, \tau_{n-1})
\]
with the same coproduct formulas as in $\dual{A}$. We may also write
$\A_{< \infty}$ for $\A$.
\item For any sub-Hopf algebra $\alg{B}$ of $\A$, we write
$\alg{B}_{\geq n}$ for $\alg{B} \cap A_{\geq n}$, and similarly
$\alg{B}_{<n} = \alg{B} \modmod \alg{B}_{\geq n}$.
\item We will not need to use elementary sub-Hopf algebras of
$\A$ in any crucial way in this paper, but we include the following
description for reference and to back up the claim in
the introduction that elementary sub-Hopf algebras have bounded
profile functions and lie in $\D$. The following facts are proved in
\cite[2.1.7--8]{palmieri;steenrod-stable}. For each $r \geq -1$ there
is a maximal elementary sub-Hopf algebra $\alg{E}(r)$ given by the
profile function
\begin{align*}
\alg{E}(-1): & \ ((0, 0, 0, \dots); (2, 2, 2, \dots)), \\
\intertext{and for $r \geq 0$}
\alg{E}(r): & \ ((\underbrace{0, \dots, 0}_{r}, r+1, r+1, r+1, \dots);
(\underbrace{1, \dots, 1}_{r+1}, 2, 2, 2, \dots))
\end{align*}
We denote the profile function of $\alg{E}(r)$ by 
\[
((f^{r}_{1}, f^{r}_{2}, \dots); (e^{r}_{0}, e^{r}_{1}, \dots)).
\]
If for some $r$, the list $((b_{1}, b_{2}, \dots); (e_{0},
e_{1}, \dots))$ satisfies $b_{i} \leq f^{r}_{i}$ and $e_{j} \leq
e^{r}_{j}$ for all $i, j$, then this is a valid profile function and
corresponds to an elementary sub-Hopf algebra of $\alg{E}(r)$.  The
maximal elementary sub-Hopf algebras are normal; the non-maximal ones
need not be.
\end{itemize}

\begin{remark}\label{remark-adams-margolis-A<n}
We note that Adams and Margolis have classified the sub-Hopf
algebras of $\A_{<n}$: see \cite[Theorem
3.3]{adams-margolis;alg}. These have the same form as for $\A$: their
duals look like $\dual{A}_{<n} / (\xi_{1}^{p^{b_{1}}}, \dots,
\xi_{n-1}^{p^{b_{n-1}}}; \tau_{0}^{e_{0}}, \dots,
\tau_{n-1}^{e_{n-1}})$ for exponents $((b_{1}, \dots,
b_{n-1}); (e_{0}, \dots, e_{n-1}))$ satisfying the same conditions as
in \autoref{thm-adams-margolis} above.
\end{remark}


\subsection{Cohomology and Steenrod operations}\label{subsec-ext}

We will usually write $H^{*}(\alg{B};-)$ for $\Ext_{\alg{B}}^{*}(\Fp,
-)$, and we may omit the coefficients if they are trivial:
$H^{*}(\alg{B})$ means $H^{*}(\alg{B};\Fp)$.

As codified by May in \cite[Theorem 11.8]{may;steenrod},
there are Steenrod operations acting on the mod $p$ cohomology of a
cocommutative Hopf algebra $\alg{B}$, and in particular for sub-Hopf
algebras of $\A$. We recall May's notation:
\begin{gather*}
\PP{i}: H^{s,t}(\alg{B}) \to H^{s+(2i-t)(p-1), pt}(\alg{B}), \\
\beta \PP{i}: H^{s,t}(\alg{B}) \to H^{s+1+(2i-t)(p-1), pt}(\alg{B}).
\end{gather*}
These satisfy the Cartan formula as well as some instability conditions:
\begin{gather*}
\PP{n}(xy) = \sum_{i+j=n} \PP{i}(x) \otimes \PP{j}(y), \\
\beta \PP{n}(xy) = \sum_{i+j=n} (\beta \PP{i}(x) \PP{j}(y) + \PP{i}(x)
\beta \PP{j}(y)), \\
\PP{i} = 0 \quad \text{if $2i <t$ or $2i > s+t$}, \\
\beta \PP{i} = 0 \quad \text{if $2i <t$ or $2i \geq s+t$}, \\
\PP{i}(x) = x^{p} \quad \text{if $2i=s+t$}.
\end{gather*}
We also recall May's regrading \cite[p.\@ 227]{may;steenrod} of these
operations:
\begin{align*}
\tildePP{i} &= \PP{i+t}: H^{s,2t} \to H^{s+2i(p-1),2pt}, \\
\beta \tildePP{i} &= \beta \PP{i+t}: H^{s,2t} \to H^{s+2i(p-1)+1,2pt}.
\end{align*}
As May points out, these regraded operations only act on $H^{s,t}$
when $t$ is even, so they are not useful when dealing with the
exterior part of the dual Steenrod algebra and any cohomology classes
arising therefrom.

When the internal degree $t$ is even, the Steenrod operation
$\PP{t/2}: H^{s,t} \to H^{s,pt}$ is analogous to the operation
$\PP{0}$ on the cohomology of a space, and indeed with the regrading,
this operation is $\tildePP{0}$. It is not, however, the identity in
general: it is instead induced by the $p$th power map on the dual ---
see \cite[Proposition 11.10]{may;steenrod}.


\subsection{The algebra \texorpdfstring{$D(x)$}{D(x)} and the elements
\texorpdfstring{$h_{ts}$, $b_{ts}$, and $a_{n}$}{hts, bts, and an}}

Despite the unfortunate notation clash with $\D$, we follow
Miller-Wilkerson \cite{miller-wilkerson} and write $D(x)$ for the
truncated polynomial algebra $D(x) = \Fp [x] / (x^{p})$, and we write
$D(x_{1}, x_{2}, \dots)$ for $D(x_{1}) \otimes D(x_{2}) \otimes \cdots
$. It is standard that when $p$ is odd, $H^{*}(D(x)) \cong
\exterior_{\Fp}(h) \otimes \Fp [b]$, where $h \in H^{1, |x|}$ and $b
\in H^{2, p|x|}$, and $b$ can be obtained from $h$ as a $p$-fold
Massey product, $b = \massey{h, \dots , h}$. If $D(x)$ is given the
structure of a cocommutative Hopf algebra, e.g.\@ with $x$ primitive,
then one has Steenrod operations acting on its cohomology, and in that
case $b$ can also be obtained via $b = \beta \PP{|h|/2}(h)$.

The Milnor basis for $\A$ is obtained by dualizing the monomial basis
for $\dual{A}$. Let $\pst{s}{t}$ be the Milnor basis element dual to
$\xi_{t}^{p^{s}}$. If $\alg{B}$ is a sub-Hopf algebra of $\A$ and if
$\xi_{t}^{p^{s}}$ is primitive in its dual $\dual{B}$, then the
algebra $D(\pst{s}{t})$ is a quotient algebra of $\alg{B}$, and so
there is an induced map
\[
H^{*}(D(\pst{s}{t})) \to H^{*}(\alg{B}).
\]
We denote the algebra generators in the domain by $h_{t,s} \in H^{1,
|\xi_{t}^{p^{s}}|}(D(\pst{s}{t}))$ and $b_{t,s} \in H^{2,
p|\xi_{t}^{p^{s}}|}(D(\pst{s}{t}))$, and we use the same names for
their images in $H^{*}(\alg{B})$. Of course algebra generators are
only defined up to a scalar multiple, and it might be useful to be
more precise, so we define $h_{t,s}$ to be cohomology class
represented by the cobar element $[\xi_{t}^{p^{s}}]$, and we define
$b_{t,s}$ by $b_{t,s} = \beta \PP{|\xi_{t}^{p^{s}}|/2}(h_{t,s}) =
\beta \tildePP{0}(h_{t,s})$. We often omit the commas, writing
$h_{ts}$ and $b_{ts}$ instead.

If $\xi_{t}^{p^{s}}$ and $\xi_{t}^{p^{s+1}}$ are both primitive in
$\dual{B}$, then because $\tildePP{0}$ is induced by the $p$th power
map on the dual, we have $\tildePP{0}(h_{ts}) = h_{t,s+1}$ and
$\tildePP{0}(b_{ts}) = b_{t,s+1}$ in the cohomology of $\alg{B}$.

We also recall that $H^{*}(\exterior (x)) \cong \Fp [a]$
where $a \in H^{1, |x|}$ is represented by $[x]$ in the cobar
complex. For each $n$, we have a sub-Hopf algebra $\exterior (Q_{n})
\subseteq \A$, dual to the quotient $\exterior (\tau_{n})$ of
$\dual{A}$. We let $a_{n}$ denote the algebra generator of
$H^{*}(\exterior (Q_{n}))$. If $\alg{B}$ is a sub-Hopf
algebra of $\A$ with $\tau_{n}$ primitive in $\dual{B}$, then
$\exterior (Q_{n})$ is a quotient Hopf algebra of $\alg{B}$, and we
also use $a_{n} \in H^{*}(\alg{B})$ to denote the image of $a_{n}
\in H^{*}(\exterior (Q_{n}))$.


\subsection{Two lemmas}

We present two lemmas which will be used throughout the paper.

As above, given $D(x) = \Fp[x]/(x^{p})$, we write $h, b$ for the
algebra generators of $H^{*}(D(x)) \cong \exterior (h) \otimes \Fp
[b]$. (This is assuming that $p$ is odd; if one wants to consider this
when $p=2$, use the $\exterior (x)$ case below.)  Given an extension
of Hopf algebras over $\Fp$
\[
\alg{C} \to \alg{B} \to D(x),
\]
there is an induced map $H^{*}(D(x)) \to H^{*}(\alg{B})$, and we also
write $h,b$ for the images of those generators under this map. 
Similarly, if we have a Hopf algebra extension
\[
\alg{C} \to \alg{B} \to \exterior (x),
\]
there is an induced map $H^{*}(\exterior (x)) \to H^{*}(\alg{B})$, and
we use $a$ to denote both a generator of $H^{*}(\exterior (x))$ and
its image in $H^{*}(\alg{B})$. The following is due to Quillen and Venkov
\cite{quillen-venkov}.

\begin{lemma}[Quillen-Venkov]\label{lemma-quillen-venkov}
\begin{enumerate}
\item Given a Hopf algebra extension $\alg{C} \to \alg{B} \to D(x)$
as above, if $z \in H^{*}(\alg{B})$ is in the kernel of the restriction
map to $H^{*}(\alg{C})$, then $z^{2}$ is in the ideal generated by $b$.
\item Given a Hopf algebra extension $\alg{C} \to \alg{B} \to \exterior (x)$
as above, if $z \in H^{*}(\alg{B})$ is in the kernel of the restriction
map to $H^{*}(\alg{C})$, then $z$ is in the ideal generated by $a$.
\end{enumerate}
\end{lemma}

\begin{proof}
Quillen and Venkov's result is about group cohomology, but their proof
carries over essentially unchanged to this setting.  In a bit more
detail, for (a) the stated setup for their lemma is a little
different: they start with a nonzero element $h \in H^{1}(G;\Z/p\Z) =
\Hom(G, \Z/p\Z)$, and then view $h$ as corresponding to a group
extension $1 \to G' \to G \to \Z/p\Z \to 1$, in which case $h$ is the
image of the fundamental class in $H^{1}(\Z/p\Z;\Z/p\Z)$ --- the
element $h$ is doing double duty here. Then they show that if $z \in
H^{*}(G;\Z/p\Z)$ restricts to zero in $H^{*}(G';\Z/p\Z)$, then $z^{2}$ is
in the ideal generated by $\beta (h)$, where $\beta$ is the Bockstein
homomorphism.

In our setup we start with a Hopf algebra extension $\alg{C} \to
\alg{B} \to D(x)$ and we \emph{choose} $h$ to be the image in
$H^{1}(\alg{B})$ of any nonzero element in $H^{1}(D(x))$.  Given this
setup, their proof works.

Their proof works for (b) also, or alternatively, their Remark on
p.\@ 317 applies. Indeed, according to \cite[Lemma
1.2.15(a)]{palmieri;steenrod-stable}, when we have an extension
$\alg{C} \to \alg{B} \to \exterior (x)$, then we get a long exact
sequence
\[
\cdots \to H^{*}(\alg{B}) \xrightarrow{\cdot a} H^{*}(\alg{B})
\xrightarrow{i^{*}} H^{*}(\alg{C}) \to \cdots,
\]
where $i^{*}$ is the restriction map. So $z$ is in the kernel of
$i^{*}$ if and only if it is in the image of multiplication by $a$.
\end{proof}

See also \cite[Theorem 5.2.1]{benson;cup2} for a different proof of
the group theory result which also translates well to this situation;
see \cite[Proposition 2.4]{lin;jpaa-I} for a version for Hopf
algebras, but only in characteristic two; and see also \cite[Lemma
3.2]{palmieri;quasi} for a closely related result.

The following is an odd primary version of a result of W. H. Lin
\cite[Theorem 3.1]{lin;jpaa-I} and our proof is modeled on Lin's. His
result is stronger: he shows that when $p=2$, the analog of $b_{ts}
b_{nk}$ is nilpotent if $k\geq t$, whereas we need to assume that $k
\geq s+t+1$. This discrepancy is at least partly unavoidable:
Wilkerson showed that at odd primes, $b_{10} b_{21}$ is non-nilpotent
(in the setting of the lemma), but there are cases which are not
covered by the lemma and where the nilpotence of $b_{ts} b_{nk}$ is
not known. We discuss this in \autoref{sec-discussion}.

\begin{lemma}[Nilpotence lemma]\label{lemma-nilpotence}
Fix $m \leq \infty$ and let $\alg{B}$ be a sub-Hopf algebra of
$\A_{<m}$ with graded dual $\dual{B}$. Fix positive integers $t$ and
$n$ with $t<n$ and $t+n<m$.
\begin{enumerate}
\item Fix an integer $s\geq 0$ and assume that
$\xi_{i}^{p^{s}}=0$ in $\dual{B}$ for $i<t$, and assume that
$\xi_{t}^{p^{s}}\neq 0$ while $\xi_{t}^{p^{s+1}} = 0$.  Fix $k \geq
s+t+1$, and assume that $\xi_{i}^{p^{k}} =0$ for $i<n$. Finally,
assume that $\xi_{n}^{p^{k}} \neq 0$. Then $b_{ts} b_{nk}$ is
nilpotent in $H^{*}(\alg{B})$.
\item Assume that
$\tau_{i}=0$ in $\dual{B}$ for $i<t$, and assume that $\tau_{t} \neq 0$.
Fix $k \geq t$, and assume that $\xi_{i}^{p^{k}} =0$ for $i<n$. Finally,
assume that $\xi_{n}^{p^{k}} \neq 0$. Then $a_{t} b_{nk}$ is
nilpotent in $H^{*}(\alg{B})$.
\end{enumerate}
\end{lemma}

Here is a sketch of the polynomial part of the profile function of
$\alg{B}$ in part (a), where terms in \textcolor{gray}{gray} are
zero in $\dual{B}$:

\medskip
\begin{tikzpicture}[scale=0.7]
\draw (0,7.5) -- (0,0) -- (10,0);
\draw (4.5, 0.5) node {$\xi_{t}$};
\draw (4.5, 1.5) node {$\vdots$};
\draw (4.5, 2.5) node {$\xi_{t}^{p^{s}}$};
\draw[gray] (3.5, 2.5) node {$\xi_{t-1}^{p^{s}}$};
\draw[gray] (4.5, 3.5) node {$\xi_{t}^{p^{s+1}}$};
\draw (4, 2) -- (4, 3) -- (5, 3);
\draw (7, 5.7) -- (7, 7);
\draw (7.5, 0.5) node {$\xi_{n}$};
\draw (7.5, 1.5) node {$\vdots$};
\draw[gray] (6.5, 6.5) node {$\xi_{n-1}^{p^{k}}$};
\draw (7.5, 6.5) node {$\xi_{n}^{p^{k}}$};
\draw (6, 4.2) node {$\iddots$};
\end{tikzpicture}

\begin{proof}
(a) Note that the assumptions guarantee that $\xi_{t}^{p^{s}}$ and
$\xi_{n}^{p^{k}}$ are primitive, so that $b_{ts}$ and $b_{nk}$ are
elements of $H^{2}(\alg{B})$.  Also, the hypothesis that $t+n<m$,
together with the fact that $\alg{B}$ is a sub-Hopf algebra of
$\A_{<m}$, guarantees that $\xi_{t+n} \neq 0 \in \dual{B}$. This is
easiest to see from the dual point of view: if $\xi_{t+n} = 0$, then
$\Delta (\xi_{t+n}) = 0$, but $\Delta (\xi_{t+n})$ contains a term
$\xi_{n}^{p^{t}} \otimes \xi_{t}$, and this term is nonzero in
$\dual{B}$. Indeed, after raising everything to the $p^{s}$ power,
this is the only term in the reduced coproduct on $\xi_{t+n}^{p^{s}}$:
we have
\[
\xi_{t+n}^{p^{s}} \mapsto \xi_{n}^{p^{s+t}} \otimes \xi_{t}^{p^{s}}.
\]
We proceed by induction, first considering the case when $k=s+t+1$.
We note that $\xi_{n}^{p^{s+t}}$ and $\xi_{t}^{p^{s}}$ are both
primitive in $\dual{B}$, so both $h_{n,s+t} = [\xi_{n}^{p^{s+t}}]$ and
$h_{ts}=[\xi_{t}^{p^{s}}]$ are cohomology classes. This is clear for
$\xi_{t}^{p^{s}}$. For $\xi_{n}^{p^{s+t}}$, its reduced coproduct is
\[
\sum_{i=t}^{n-t} \xi_{n-i}^{p^{s+t+i}} \otimes \xi_{i}^{p^{s+t}}.
\]
We are assuming that $\xi_{j}^{p^{k}} = \xi_{j}^{p^{s+t+1}}=0$ for all
$j<n$, so the first tensor factor in each summand is zero. Thus the
coproduct on $\xi_{t+n}^{p^{s}}$ produces a relation in cohomology:
\[
h_{n,s+t} h_{ts} = 0.
\]
Applying Steenrod operations to this relation gives the following ---
we label each line with the operation being applied:
\begin{align*}
\beta \tildePP{0}: & \quad b_{n,s+t} h_{t,s+1} - h_{n,s+t+1} b_{ts} = 0, \\
\intertext{Note that the first term is zero since $\xi_{t}^{p^{s+1}}=0$.}
\beta \tildePP{1}: & \quad - b_{n,s+t+1} b_{ts}^{p} = 0.
\end{align*}
This finishes the case when $k=s+t+1$.

The same argument does not immediately apply if $k>s+t+1$: we would
again want to consider the coproduct on $\xi_{t+n}^{p^{s}}$, but the
resulting tensor factors need not be primitive and so would not
correspond to cohomology classes, so we would not be able to apply
Steenrod operations to them. So we first work in a smaller Hopf
algebra in which the terms are primitive, and then we enlarge the Hopf
algebra. 

The inductive hypothesis is the following: for each $k'$ with $s+t+1
\leq k' \leq k$, in the cohomology of any Hopf algebra whose dual
satisfies
\begin{enumerate}
\item $\xi_{i}^{p^{s}}=0$ for $i<t$,
\item $\xi_{t}^{p^{s}} \neq 0$, $\xi_{t}^{p^{s+1}}=0$,
\item $\xi_{n}^{p^{k}} \neq 0$,
\item $\xi_{i}^{p^{k'}}=0$ for $t<i<n$,
\end{enumerate}
the element $b_{t0}b_{n,k'}$ is nilpotent. We have verified this when
$k=s+t+1$, starting the induction. Note also that the argument that we
give is independent of $n$, so when we apply the inductive hypothesis,
we can (and will) apply it with different values of $n$.

For the inductive step we first work in a Hopf algebra $\alg{C}$ whose
dual satisfies (a), (b), and (c), together with
\begin{enumerate}
\item [(d$'$)] $\xi_{i}^{p^{s+t}} \neq 0$, $\xi_{i}^{p^{s+t+1}} = 0$ for $t<i<n$.
\end{enumerate}
We will show that $b_{ts} b_{nk}$ is nilpotent in $H^{*}(\alg{C})$;
after that, we will have to modify $\alg{C}$ to get a Hopf algebra
that sets us up for the next step in the induction.  Regarding
hypothesis (d$'$) and in particular the requirement that
$\xi_{i}^{p^{s+t}} \neq 0$: any Hopf algebra $\alg{C}'$ whose dual
satisfies (a), (b), (c), and $\xi_{i}^{p^{s+t+1}}=0$ for $t<i<n$ will be a
sub-Hopf algebra of $\alg{C}$; thus as long as $b_{ts}b_{nk}$ is
nilpotent in $H^{*}(\alg{C})$, it will also be nilpotent in
$H^{*}(\alg{C}')$ via the restriction map.

The reduced coproduct in $\alg{C}$ on $\xi_{t+n}^{p^{s}}$ is
\[
\sum_{i=1}^{n+t-i} \xi_{n+t-i}^{p^{s+i}} \otimes \xi_{i}^{p^{s}}.
\]
Since $\xi_{i}^{p^{s}}=0$ for $i<t$, the lower limit on the sum can be
changed to $t$. Similarly, since $\xi_{n+t-i}^{p^{s}}=0$ when
$n+t-i<t$, then the upper limit can be changed to $n$: the reduced
coproduct is
\[
\sum_{i=t}^{n} \xi_{n+t-i}^{p^{s+i}} \otimes \xi_{i}^{p^{s}}.
\]
Since $\xi_{j}^{p^{s+t+1}}=0$ when $t<j<n$, the $i=t$ term will be the
only one present: the coproduct is $\xi_{n}^{p^{s+t}} \otimes
\xi_{t}^{p^{s}}$. Both tensor factors here are primitive, so we get a
relation in cohomology: $h_{n,s+t}h_{ts}=0$. We can then apply
Steenrod operations as above to get the desired result:
\begin{align*}
\beta \tildePP{0}: & \quad b_{n,s+t} h_{t,s+1} - h_{n,s+t+1} b_{ts} = 0, \\
\intertext{Note that the first term is zero since $\xi_{t}^{p^{s+1}}=0$.}
\beta \tildePP{1}: & \quad - b_{n,s+t+1} b_{ts}^{p} = 0, \\
\tildePP{p}: & \quad -b_{n,s+t+2} b_{ts}^{p^{2}} = 0, \\
\vdots \\
\tildePP{p^{j-1}}: & \quad -b_{n,s+t+j} b_{ts}^{p^{j}} = 0.
\end{align*}
So by choosing an appropriate $j$, we see that $b_{ts} b_{nk}$ is
nilpotent in $H^{*}(\alg{C})$.

Now we need to get a nilpotence result over a different Hopf algebra,
one whose dual satisfies (a), (b), and (c) above, along with
\begin{itemize}
\item[(d$''$)] $\xi_{i}^{p^{k}} = 0$ for $t<i<n$.
\end{itemize}
As noted above, we may assume that $\xi_{i}^{p^{s+t}} \neq 0$ for
$t<i<n$, and we can build such a Hopf algebra from $\alg{C}$ by adding
in one ``box'' $\xi_{i}^{p^{j}}$ at a time: we assume that the result
holds over a Hopf algebra $\alg{C}'$ and we want to show that it
continues to hold over $\alg{C}''$, where there is a Hopf algebra
extension
\[
\alg{C}' \to \alg{C}'' \to D(\pst{j}{i}),
\]
where $t<i<n$ and $s+t+1 \leq j < k$. (Recall that $\pst{j}{i}$ is dual
to $\xi_{i}^{p^{j}}$.) We use the Quillen-Venkov lemma and our
inductive hypothesis: given $z = b_{ts} b_{nk} \in H^{*}(\alg{C}'')$,
since $z^{M} \mapsto 0$ in $H^{*}(\alg{C}')$, then $z^{2M}$ is in the
ideal generated by $b_{ij}$: we have
\[
b_{ts}^{2M} b_{nk}^{2M} = b_{ij} y
\]
for some $y$. Now because $s+t+1 \leq j < k$, by induction (using a
lower value for $n$) we know that $b_{ts} b_{ij}$ is nilpotent, so
after multiplying both sides by $b_{ts}$ and raising to a large enough
power, we see that the right-hand side is nilpotent, and hence so is
the left. This completes the inductive step.

(b) This proof is very similar. First assume that $\xi_{i}^{p^{t}} =
0$ when $i<n$. Then we have the reduced coproduct $\tau_{t+n} \mapsto
\xi_{n}^{p^{t}} \otimes \tau_{t}$, which gives us the relation $h_{nt}
a_{t} = 0$ in cohomology. If we apply the operation $\beta
\PP{p^{t}}$, we get $b_{nt} a_{t}^{p}=0$, and then applying further
operations, we get $b_{n,t+i} a_{t}^{p^{i+1}}=0$ for all $i \geq 0$.

Now we use the Quillen-Venkov lemma to get from the hypothesis that
$\xi_{i}^{p^{t}}=0$ to $\xi_{i}^{p^{k}}=0$ for $i<n$. As in part (a),
we will have Hopf algebra extensions of the form
\[
\alg{C}' \to \alg{C}'' \to D(\pst{j}{i})
\]
where $i<n$ and $t\leq j<k$, and we assume that $z = b_{nk} a_{t} \in
H^{*}(\alg{C}'')$ maps to a nilpotent element in $H^{*}(\alg{C}')$,
say $z^{M} \mapsto 0$. Then in $H^{*}(\alg{C}'')$, the Quillen-Venkov
lemma tells us that $z^{2M} = b_{nk}^{2M} a_{t}^{2M} = b_{ij}y$ for
some $y$. By induction we know that $b_{ij} a_{t}$ is nilpotent, and
that finishes the proof.
\end{proof}

Part (a) of this lemma can be improved. We have stated the version
that we need, but it is not too hard to get the same conclusion with
slightly weaker hypotheses:

\begin{lemma}\label{lemma-nilpotence-stronger}
Fix $m \leq \infty$ and let $\alg{B}$ be a sub-Hopf algebra of
$\A_{<m}$ with graded dual $\dual{B}$. Fix positive integers $t$ and
$n$ with $t<n$ and $t+n<m$. Fix an integer $s\geq 0$.  Assume that
$\xi_{t}^{p^{s}}\neq 0$ while $\xi_{t}^{p^{s+1}} = 0$. Fix $k \geq
s+t+1$, and assume that $\xi_{i}^{p^{k}} =0$ for $i<n$. Assume that
for $i<t$, $\xi_{i}^{p^{k-i}}=0$. Assume that $\xi_{t}^{p^{s}}$ is
primitive in $\alg{B}$. Finally, assume that
$\xi_{n}^{p^{k}} \neq 0$. Then $b_{ts} b_{nk}$ is nilpotent in
$H^{*}(\alg{B})$.
\end{lemma}

The difference is that for $i<t$, in the Nilpotence Lemma we assume
that $\xi_{i}^{p^{s}}=0$, while here we assume that
$\xi_{i}^{p^{k-i}}=0$. Since $k\geq s+t+1 > s+i+1$, then $k-i \geq s$,
so this is a weaker assumption. The original assumption that
$\xi_{i}^{p^{s}}=0$ for $i<t$ implies that $\xi_{t}^{p^{s}}$ is
primitive, while that may not follow from this weaker assumption, so
we add the primitivity assumption explicitly here. Perhaps the lemma
could be improved further: it would be nice if one only had to assume
that $\xi_{t}^{p^{s}}$ and $\xi_{n}^{p^{k}}$ were primitive.

\begin{proof}[Sketch of proof] We use the Quillen-Venkov Lemma and the
Nilpotence Lemma. Suppose that $\alg{B}$ is as in the statement, and
let $(b_{1}, b_{2}, \dots )$ be its profile function. We define a
sub-Hopf algebra $\alg{C}$ by setting $c_{i} = \min (b_{i}, s)$ if
$i<t$, $c_{i}=b_{i}$ if $i \geq t$. Then $\alg{C}$ satisfies the
hypotheses of the Nilpotence Lemma, so $b_{ts} b_{nk}$ is nilpotence
in the cohomology of $\alg{C}$. We can then build $\alg{B}$ from
$\alg{C}$ by using induction and adding in one element $\xi_{i}^{p^{j}}$
at a time. That is, as in the above proof, in the inductive step we
have a Hopf algebra extension
\[
\alg{C}' \to \alg{C}'' \to D(\pst{j}{i}),
\]
where $i<t$, $j<k-i$, so $i+j+1 \leq k$. We assume that $b_{ts}
b_{nk}$ is nilpotence over $\alg{C}'$, so by the Quillen-Venkov Lemma,
some power of $b_{ts} b_{nk}$ is a multiple of $b_{ij}$. Since $k \geq
i+j+1$, $b_{ij} b_{nk}$ is nilpotent, so $b_{ts} b_{nk}$ remains
nilpotent over $\alg{C}''$.
\end{proof}


\section{Reducing to nilpotence of particular cohomology classes}
\label{sec-reduce}

The goal of this section is to prove \autoref{thm-reduction}, the
point of which is to reduce the questions of detection of nilpotence
and existence of an $F$-isomorphism to specific cohomology
calculations: are certain cohomology classes nilpotent? We need a bit
of setup and notation first, so suppose that we have Hopf algebras
$\alg{C} \subseteq \alg{H} \subseteq \A$. Suppose further that
$\alg{C}$ and $\alg{H}$ agree on the exterior part of the Steenrod
algebra: we may write their profile functions as
\begin{align*}
\alg{C} &:
((c_{1}, c_{2}, c_{3}, \dots); 
(e_{0}, e_{1}, e_{2}, \dots))\\
\alg{H} &:
((h_{1}, h_{2}, h_{3}, \dots);
(e_{0}, e_{1}, e_{2}, \dots)),
\end{align*}
with $c_{j} \leq h_{j}$ for all $j$. Given $n$, we write $\alg{C}_{n}$
for the sub-Hopf algebra of $\alg{H}$ with profile function
\[
((c_{1}, \dots, c_{n-1}, h_{n}, h_{n+1}, \dots); 
(e_{0}, e_{1}, e_{2}, \dots))
\]
(Note that $\alg{C}_{n}$ is not perfect notation since the algebra
depends on $\alg{H}$ in addition to $\alg{C}$ and $n$, but we think of
$\alg{H}$ as being the ambient Hopf algebra and suppress it from the
notation.)  Here is a picture of the polynomial part of the profile
functions:

\medskip
\begin{tikzpicture}[scale=0.7]
\draw (0,7.5) -- (0,0) -- (10,0);
\draw (0,3) -- (9, 7);
\draw[above] (4, 4.8) node {$\alg{H}$};
\draw (0,0) -- (9, 4);
\draw[below] (7, 3) node {$\alg{C}$};
\draw [thick] (0,0.1) -- (5,2.33) -- (5,5.12) -- (9,6.9);
\draw [right] (5, 4.2) node {$\alg{C}_{n}$};
\draw [below] (5.5, 0) node {$n$};
\draw [dotted] (5,0) -- (5,2.3);
\draw (5.3, 2.7) node {$\bullet$};
\draw (8,5) node {$\xi_{n}^{p^{c_{n}}}$};
\draw [->, shorten <=5mm, shorten >= 2mm] (8,5) -- (5.3, 2.7);
\end{tikzpicture}

For example, we will apply this in the next section with $\alg{H} =
\A$ and $\alg{C} = \Binom$, in which case the algebra $\Binom_{n}$ has
profile function
\[
((1, 2, 4, \dots,
\text{\footnotesize{$\binom{n-1}{2}+1$}}, 
\infty, \infty, \infty, \dots);
(2, 2, 2, \dots)).
\]

\begin{theorem}\label{thm-reduction}
Fix $n \geq 1$. Suppose that we have Hopf algebras $\alg{C} \subseteq
\alg{H} \subseteq \A$, with $\alg{C}$ agreeing with $\alg{H}$ on the
exterior part of $\A$. Suppose that $\alg{C}$ is normal in $\alg{H}$,
and define $\alg{C}_{n}$ as above. Write the profile functions for
$\alg{C}$ and $\alg{H}$ as $((c_{1}, c_{2}, \dots); (e_{0}, e_{1},
\dots))$ and $((h_{1}, h_{2}, \dots); (e_{0}, e_{1}, \dots)$,
respectively. Suppose that for each $n$, either $c_{n}=h_{n}$ or
$b_{n,c_{n}}$ is nilpotent in the cohomology of $\alg{C}_{n}$. Then
\begin{enumerate}
\item $\alg{C}$ detects nilpotence in $\stable{H}$.
\item The map $H^{*}(\alg{H}) \to \left( H^{*}(\alg{C})
\right)^{\alg{H}}$ is an $F$-isomorphism.
\end{enumerate}
\end{theorem}

(We could write (b) in terms of $(H^{*}(\alg{C}))^{\alg{H} \modmod
\alg{C}}$, but of course that is equal to
$(H^{*}(\alg{C}))^{\alg{H}}$. Here and below we may use the shorter
version for more compact notation.)

Note that if $c_{n} \neq h_{n}$, then $c_{n} < \infty$. In this
situation, $\xi_{n}^{p^{c_{n}}}$ will be nonzero and primitive in
$\dual{C}_{n}$, so $b_{n,c_{n}}$ will be a cohomology class. (Why is
it primitive?  $\xi_{n}^{p^{c_{n}}}$ is zero in $\dual{C}$, so its
reduced coproduct is zero in $\dual{C} \otimes \dual{C}$. That
coproduct is still zero in $\dual{C}_{n} \otimes \dual{C}_{n}$, which
means that $\xi_{n}^{p^{c_{n}}}$ is primitive there.)

We also point out that if $b_{n,c_{n}}$ is nilpotent, then so is
$b_{n,j}$ for any $j \geq c_{n}$, because the Steenrod operation
$\tildePP{0}$ is an algebra map sending $b_{n,i}$ to $b_{n,i+1}$ for
all $i$.

The proof of the theorem is almost identical to the proofs of
\cite[Theorems 5.1.2 and 5.1.5]{palmieri;steenrod-stable} when $p=2$.
We assume that $c_{n} < h_{n}$, so $\xi_{n}^{p^{c_{n}}} \neq 0$ in
$\dual{C}_{n}$. For each $n$, $\alg{C}_{n+1}$ is normal in
$\alg{C}_{n}$, and we have
\[
\alg{C}_{n} \modmod \alg{C}_{n+1} \cong D(\pst{c_{n}}{n},
\pst{c_{n}+1}{n}, \dots).
\]
This has the structure of a $\alg{C}_{n}$-module, and for any $j \geq
c_{n}$, we have a finite-dimensional quotient $\alg{C}_{n}$-module
\[
W_{j} = D(\pst{c_{n}}{n}, \pst{c_{n}+1}{n}, \dots, \pst{j}{n}).
\]
It also makes sense to let $W_{c_{n}-1} = \Fp$: in the proof of the
following lemma we will use the observation that each $W_{j+1}$ can be
built from $W_{j}$, and this process starts with $W_{c_{n}-1} =
\Fp$. We abuse notation and also use the same notation $W_{j}$ for the
corresponding object in $\cat{Stable}(\dual{C}_{n})$.

\begin{lemma}\label{lemma-thick}
We work in the category $\cat{Stable}(\dual{C}_{n})$.  Assume that for
some $j \geq c_{n}$, the element $b_{nj} \in H^{*}(\alg{C}_{n})$ is
nilpotent. Then the object $W_{j-1}$ is in the thick subcategory
generated by $W_{j}$. So if $b_{nj}$ is nilpotent for all $j \geq
c_{n}$, then $\Fp$ is in the thick subcategory generated by any
$W_{j}$.
\end{lemma}

This lemma is where the nilpotence of $b_{n,c_{n}}$ is used.

\begin{proof}
By \cite[Lemma 1.2.15(b)]{palmieri;steenrod-stable}, for any $j \geq
c_{n}$ there are two cofiber sequences in
$\cat{Stable}(\dual{C}_{n})$:
\begin{gather*}
\Sigma^{2, p|\xi_{n}^{p^{j+1}}|} W_{j-1} \xrightarrow{b_{nj}} W_{j-1} \to
\widetilde{W}_{j}, \\
\Sigma^{1, |\xi_{n}^{p^{j+1}}|} W_{j} \to \widetilde{W}_{j} \to W_{j}.
\end{gather*}
Therefore $\widetilde{W}_{j}$ is in the thick subcategory generated
by $W_{j}$, and the nilpotence of $b_{nj}$ tells us that $W_{j-1}$ is
in the thick subcategory generated by $\widetilde{W}_{j}$.
\end{proof}

\begin{lemma}\label{lemma-vanish}
Fix $j \geq 0$ and let $\alg{B} = D(\pst{j}{n}, \pst{j+1}{n},
\pst{j+2}{n}, \dots) = \plaindual{(\Fp [\xi_{n}^{p^{j}}])}$. If $M$ is
a $\alg{B}$-module which is zero below degree $t$, then
$H^{*,*}(\alg{B}; M)$ has a vanishing plane: there is a number $c$,
independent of $M$, so that $H^{s,u}(\alg{B};M) = 0$ when $u <
\frac{p}{2}|\xi_{n}^{p^{j}}| s + t + c$.
\end{lemma}

\begin{proof}
This is an easy computation. We know that $H^{*}(\alg{B}; \Fp) \cong
\exterior (h_{ni} : i \geq j) \otimes \Fp [b_{ni} : i \geq j]$
has a vanishing plane with a slope given by the element $b_{nj}$,
namely $m_{j} := \frac{p}{2} |\xi_{n}^{p^{j}}|$. That is,
$H^{s,u}(\alg{B}; \Fp) = 0$ when $u<m_{j} s + c$ for some intercept
$c$.  Therefore we can construct an injective resolution of $\Fp$, the
connectivity of whose terms increases with this slope. Tensor this
with $M$ to get an injective resolution of $M$; its terms will have
connectivity increasing with the desired slope, and with intercept
$t+c$.
\end{proof}

\begin{proof}[Proof of \autoref{thm-reduction}]
We will omit some details here; see \cite[Subsections
5.2.1--2]{palmieri;steenrod-stable} for full details.

We think of building $\alg{H}$ from $\alg{C}$ using inclusions of Hopf
algebras
\[
\alg{C} \subseteq \dots \subseteq \alg{C}_{n+1} \subseteq \alg{C}_{n} \subseteq
\dots \subseteq \alg{C}_{2} \subseteq \alg{C}_{1} = \alg{H}.
\]

For the nilpotence part of the proof, we use standard techniques
(omitted here, but see \cite{palmieri;steenrod-stable}) to reduce to
the case of a connective ring spectrum $R \in \stable{H}$ and an
element $\alpha \in H^{*}(\alg{H}; R)$ which restricts to zero in
$H^{*}(\alg{C}; R)$. Since $\alg{C}$ and $\alg{C}_{n}$ agree through a
range of dimensions increasing with $n$, so do their cohomology
algebras, so $\alpha$ will also restrict to zero in
$H^{*}(\alg{C}_{n+1};R)$ for $n \gg 0$. This starts an inductive
process; the inductive step is to show that if $\alpha$ restricts to
zero over $H^{*}(\alg{C}_{n+1};R)$, then it is nilpotent in
$H^{*}(\alg{C}_{n};R)$.

Note that if $c_{n}=h_{n}$, then $\alg{C}_{n} = \alg{C}_{n+1}$, so we
only have to consider the case when $c_{n}<h_{n}$, which means that
$b_{n,c_{n}}$ is nilpotent in $H^{*}(\alg{C}_{n})$.

Given the connective ring spectrum $R$ and the Hopf algebra extension
\[
\alg{C}_{n+1} \to \alg{C}_{n} \to \alg{C}_{n} \modmod \alg{C}_{n+1},
\]
the Cartan-Eilenberg spectral sequence has the form
\[
E_{2}^{s,t,*}(R) \cong H^{s,*}(\alg{C}_{n} \modmod \alg{C}_{n+1};
H^{t,*}(\alg{C}_{n+1}; R)) \Rightarrow H^{s+t,*}(\alg{C}_{n}; R).
\]
By \cite[Proposition 1.4.9]{palmieri;steenrod-stable}, this may be
viewed as an Adams spectral sequence, and \cite[Proposition
5.2.2]{palmieri;steenrod-stable} will let us get as nice a vanishing
plane as we want at some page of this spectral sequence. We will get
into the details, but first note that there are two cases here
depending on whether any power of $\xi_{n}$ is zero in the dual of the
ambient Hopf algebra $\alg{H}$: whether $h_{n}<\infty$ or not. The
case when $h_{n}<\infty$ is easier, and we handle it later: for now we
assume that $h_{n}=\infty$.

To get the vanishing plane when $h_{n}=\infty$: for any $j \geq
c_{n}$, since $\alg{C}_{n+1}$ acts trivially on $W_{j}$, we have
\begin{align*}
E_{2}^{s,t,*}(R \Smash W_{j}) &\cong
H^{s,*}(\alg{C}_{n} \modmod \alg{C}_{n+1}; H^{t,*}(\alg{C}_{n+1}; R
\Smash W_{j})) \\
&\cong 
H^{s,*}(\alg{C}_{n} \modmod \alg{C}_{n+1}; W_{j} \Smash
H^{t,*}(\alg{C}_{n+1}; R)) \\
&\cong 
H^{s,*}(\plaindual{(\Fp [\xi_{n}^{p^{j+1}}])}; H^{t,*}(\alg{C}_{n+1}; R))
\end{align*}
where we have used a change-of-rings isomorphism at the end.  By
\autoref{lemma-vanish} above, with coefficients in a bounded below
module, this cohomology has a vanishing plane of the form
\[
E_{2}^{s,t,u}(R \Smash W_{j})=0 \quad \text{when} \quad m_{j+1} s + t - u > c
\]
for some $c$, where as above, $m_{j+1} = \frac{p}{2}
|\xi_{n}^{p^{j+1}}|$. The lemma is stated for modules, but the
existence of such a vanishing plane extends to suitably connective
objects in $\stable{H}$, so since $R$ is a connective ring spectrum,
we find that $E_{2}^{s,t,u}(R \Smash W_{j})$ has a good vanishing
plane: it will be zero when $m_{j+1} s + t - u > c$, for some
$c$. Hence by \autoref{lemma-thick} and genericity of vanishing
planes in the Adams spectral sequence --- see
\cite{hopkins-palmieri-smith} --- for some $c'$ and some $r$, we see
that $E_{r}^{s,t,u}(R)$ will be zero when $m_{j+1} s + t - u > c'$.

If $\alpha \in H^{*}(\alg{C}_{n};R)$ restricts to zero in
$H^{*}(\alg{C}_{n+1}; R)$, it will not lie on the $s=0$ edge of the
spectral sequence: it will be represented by a class $a \in
E_{2}^{s,t,u}$ with $s>0$. Now if we choose a large enough $j$ so as
to obtain a low enough vanishing plane, some power of $a$ will lie
above the vanishing plane and hence will be zero; therefore $\alpha$
will be nilpotent in the abutment, as desired. (See
\cite{palmieri;steenrod-stable} for more details.)

In the case when $h_{n}<\infty$, then 
\[
\alg{C}_{n} \modmod \alg{C}_{n+1} \cong D(\pst{c_{n}}{n},
\pst{c_{n}+1}{n}, \dots, \pst{h_{n}-1}{n}) = W_{h_{n}-1}.
\]
The spectral sequence has 
\[
E_{2}^{s,t,*}(R \Smash W_{h_{n}-1}) \cong H^{s,*}(\alg{C}_{n} \modmod
\alg{C}_{n+1}; W_{h_{n-1}} \Smash H^{t,*}(\alg{C}_{n+1}; R)).
\]
The coefficient module is free over $\alg{C}_{n} \modmod
\alg{C}_{n+1}$, so $E_{2}^{s,t,*}(R \Smash W_{h_{n}-1}) = 0$ when
$s>0$. This is a very strong vanishing plane! By \autoref{lemma-thick}
and the genericity of vanishing planes, for some $c'$ and some $r$, we
see that $E_{r}^{s,t,u}(R)$ will be zero when $s > c'$. The rest of
the proof proceeds as above. (Alternatively, one could perhaps state
and prove a version of the Quillen-Venkov lemma with coefficients in a
ring spectrum in $\stable{C}$, and then deal with one ``box''
$\xi_{n}^{p^{j}}$ at a time.)

For the $F$-isomorphism part of the proof, it remains to show that if
$y$ is in the codomain of the restriction map $\overline{\imath}^{*}$,
then some $p$th power of $y$ is in the image.  It is easy to see that
if $y \in H^{*}(\alg{C})$, then $y$ is in the image of the restriction
map $H^{*}(\alg{C}_{n}) \to H^{*}(\alg{C})$ for $n \gg 0$, so again,
the issue is going from $\alg{C}_{n+1}$ to $\alg{C}_{n}$: if $y$ is in
the invariant subalgebra $(H^{*}(\alg{C}_{n+1}))^{\alg{C}_{n}}$, then
we need to show that some power of $y$ is in the image of the
restriction map $\overline{\imath}^{*}: H^{*}(\alg{C}_{n}) \to
(H^{*}(\alg{C}_{n+1}))^{\alg{C}_{n}}$. We use the same spectral
sequence. The first thing to note is that
\[
E_{2}^{0,t,u}(\Fp) \cong H^{0,*}(\alg{C}_{n} \modmod \alg{C}_{n+1};
H^{t,*}(\alg{C}_{n+1}))
\cong (H^{t,*}(\alg{C}_{n+1}))^{\alg{C}_{n} \modmod \alg{C}_{n+1}}.
\]
This is where the invariants come into play.  Because the element $y$
is in the image of the restriction, it is represented by a class in
$E_{2}^{0,t,u}$, and using Steenrod operations acting on the spectral
sequence, one can identify the possible locations of the images of $y$
and its $p$th powers under the differentials --- see \cite[Proposition
1.4.11]{palmieri;steenrod-stable}, which says that $y^{p^{n}}$
survives to $E_{p^{n}+1}$. The possible images of the differentials of
the classes $y^{p^{n}}$ lie on a line with a positive slope in the
spectral sequence. That line will eventually go above a low enough
vanishing plane, so we see that some $y^{p^{n}}$ will be a permanent
cycle.
\end{proof}

Since $\P$ is a sub-Hopf algebra of $\A$, we can of course apply the
previous result to sub-Hopf algebras of $\P$. In fact it suffices to
work with sub-Hopf algebras of $\P$, as we now demonstrate. Fix a
sub-Hopf algebra $\alg{C} \subseteq \A$.  Just as we can view $\P$ as
both a sub- and quotient Hopf algebra of $\A$, we can view $\alg{C}
\cap \P$ as both a sub- and quotient Hopf algebra of $\alg{C}$.
%

\begin{proposition}\label{prop-P-to-A}
Given a sub-Hopf algebra $\alg{C} \subseteq \A$, suppose that
$\xi_{n}^{p^{k}}$ is primitive in $\dual{C}$, and hence in
$\plaindual{(\alg{C} \cap \P)}$. If $b_{nk} \in H^{*}(\alg{C} \cap \P)$
is nilpotent, so is $b_{nk} \in H^{*}(\alg{C})$.
\end{proposition}

\begin{proof}
The quotient map $\alg{C} \to \alg{C} \cap \P$ induces an algebra map
$H^{*}(\alg{C} \cap \P) \to H^{*}(\alg{C})$, sending $b_{nk}$ to $b_{nk}$.
\end{proof}

So given $\alg{C} \subseteq \alg{H} \subseteq \A$ as in
\autoref{thm-reduction}, it suffices to check the nilpotence of
each $b_{n,c_{n}} \in H^{*}((\alg{C} \cap \P)_{n}) = H^{*}(\alg{C}_{n}
\cap \P)$. If we can do that, we may conclude that
\begin{itemize}
\item $\alg{C}$ detects nilpotence in $\stable{H}$,
\item $\alg{C} \cap \P$ detects nilpotence in
$\cat{Stable}(\plaindual{(\alg{H} \cap \P)})$,
\item the map $H^{*}(\alg{H}) \to H^{*}(\alg{C})^{\alg{H}}$ is an
$F$-isomorphism, and
\item the map $H^{*}(\alg{H} \cap \P) \to H^{*}(\alg{C} \cap
\P)^{\alg{H} \cap \P}$ is an $F$-isomorphism.
\end{itemize}


\section{Restricting to the algebra \texorpdfstring{$\Binom$}{Binom}}\label{sec-binom}

To prove \cref{thm-nilpotence,thm-f-iso}, we will
apply \autoref{thm-reduction} twice, to each of the inclusions $\D
\subseteq \Binom$ and $\Binom \subseteq \A$, or more generally given
$\alg{C} \subseteq \A$, we apply the theorem to $\alg{C} \cap \D
\subseteq \alg{C} \cap \Binom$ and $\alg{C} \cap \Binom \subseteq
\alg{C}$, to get $F$-isomorphisms
\[
H^{*}(\alg{C}) \to \left( H^{*}(\alg{C} \cap \Binom) \right)^{\alg{C}} \to \left( \left(
H^{*}(\alg{C} \cap \D) \right)^{\alg{C} \cap \Binom} \right)^{\alg{C}} \cong \left(
H^{*}(\alg{C} \cap \D) \right)^{\alg{C}},
\]
and similarly for detection of nilpotence. This means that we have to
study the nilpotence of classes $b_{nj}$ in the cohomology algebras of
certain Hopf algebras. We focus on the inclusion $\alg{C} \cap \Binom
\subseteq \alg{C}$ in this section; in the next section we will
address $\alg{C} \cap \D \subseteq \alg{C} \cap \Binom$.

\autoref{prop-P-to-A} lets us focus on sub-Hopf
algebras of $\P$ rather than $\A$. The key to making this work is
the following, an unpublished result of Burklund, Levy, and
Pstr\k{a}gowski. Recall the notation $\P_{\geq n}$ from
\autoref{subsec-hopf}: this is the Hopf algebra with graded
dual $\dual{P}_{\geq n} = \Fp [\xi_{n}, \xi_{n+1}, \dots]$.

\begin{theorem}[\cite{burklund-levy-pstragowski}]\label{thm-blp}
The element $b_{nn} \in H^{2}(\P_{\geq n})$ is nilpotent.
\end{theorem}

We give an independent proof of this when $p=3$ --- see
\autoref{thm-3} in the appendix --- and we prove a related result
in \autoref{thm-bnn-finite-profile}.

\begin{corollary}
Suppose that $\alg{C}$ is a sub-Hopf algebra of either $\A$ or
$\P$. If $j \geq n$ and $\xi_{n}^{p^{j}} \neq 0$ in $\dual{C}$, then
the element $b_{nj} \in H^{2}(\alg{C}_{\geq n})$ is nilpotent.
\end{corollary}

\begin{proof}
If $\alg{C}$ is a sub-Hopf algebra of $\P$, then we have the
restriction map $H^{*}(\P_{\geq n}) \to H^{*}(\alg{C}_{\geq n})$, and
so $b_{nn} \in H^{*}(\alg{C}_{\geq n})$ is nilpotent. The Steenrod
operation $\tildePP{0}$ is an algebra map and sends $b_{ij}$ to
$b_{i,j+1}$ as long as these classes are defined, so we conclude that
$b_{nj} \in H^{2}(\alg{C}_{\geq n})$ is nilpotent for all $j \geq n$.

If $\alg{C}$ is a sub-Hopf algebra of $\A$, we note that $\P_{\geq n}$
is a quotient of $\A_{\geq n}$, so we get an induced map
$H^{*}(\P_{\geq n}) \to H^{*}(\A_{\geq n})$, whence $b_{nn} \in
H^{*}(\A_{\geq n})$ is nilpotent. Now as in the previous paragraph,
use Steenrod operations and the restriction map
$H^{*}(\A_{\geq n}) \to H^{*}(\alg{C}_{\geq n})$.
\end{proof}

For the rest of this section we work within the Hopf algebra $\P$, so
all profile functions have only a polynomial component. If $\alg{C}$
has profile function $(c_{1}, c_{2}, \dots)$, then $\alg{C} \cap
\Binom$ will have profile function with $n$th term $\min (c_{n},
\binom{n}{2}+1)$. Recall from the start of \autoref{sec-reduce}
that, given the inclusion $\alg{C} \cap \Binom \subseteq \alg{C}$, we
write $(\alg{C} \cap \Binom)_{n}$ for the algebra whose profile
function is the same as that for $\alg{C} \cap \Binom$ for the first
$n-1$ terms, the same as that for $\alg{C}$ afterwards:
\[
(\min (c_{1}, 1), \min (c_{2}, 2), \min (c_{3}, 4),
\dots, \min (c_{n-1}, \text{\footnotesize{$\binom{n-1}{2}+1$}}), 
c_{n}, c_{n+1}, \dots).
\]
According to \autoref{prop-P-to-A} and
\autoref{thm-reduction}, we want to prove the following.

\begin{proposition}\label{prop-bnn-binom-P}
Fix a sub-Hopf algebra $\alg{C} \subseteq \P$.  If $k \geq
\binom{n}{2}+1$ and $\xi_{n}^{p^{k}} \neq 0$ in $\dual{C}$,
then the element $b_{nk} \in H^{2}((\alg{C} \cap \Binom)_{n})$ is
nilpotent.
\end{proposition}

\begin{corollary}\label{cor-binom}
For any $\alg{C} \subseteq \A$, $\alg{C} \cap \Binom$
detects nilpotence in $\stable{C}$, and
the restriction map
\[
H^{*}(\alg{C}) \to (H^{*}(\alg{C} \cap \Binom))^{\alg{C}}
\]
is an $F$-isomorphism.
\end{corollary}

\begin{proof}[Proof of \autoref{prop-bnn-binom-P}] 
Let $\alg{B} = (\alg{C} \cap \Binom)_{n}$. Let $(b_{1}, b_{2}, \dots)$
be the profile function for $\alg{B}$, which means that $b_{i} = \min
(c_{i}, \binom{i}{2}+1)$ when $i < n$, while $b_{i} = c_{i}$ when $i
\geq n$. Note that $\alg{B}_{\geq n} = \alg{C}_{\geq n}$. Fix $k \geq
\binom{n}{2}+1$; we may assume that $\xi_{n}^{p^{k}} \neq 0$ in
$\dual{B}$ and hence also in $\dual{C}$. We use downward induction via
a sequence of sub-Hopf algebras of $\alg{B}_{n}$, first of the form
\[
\alg{C}_{\geq n} =
\alg{B}_{\geq n} \hookrightarrow
\alg{B}_{\geq n-1} \hookrightarrow \dots \hookrightarrow
\alg{B}_{\geq 2} \hookrightarrow
\alg{B}_{\geq 1} = \alg{B} = (\alg{C} \cap \Binom)_{n},
\]
and then we will interpolate further, showing that $b_{nk}$ is
nilpotent over each of them. The induction starts with
\autoref{thm-blp} and its corollary: $b_{nk}$ is nilpotent in
$H^{*}(\alg{B}_{\geq n})$. We want to show that for each $m \leq n-1$,
if $b_{nk}$ is nilpotent in $H^{*}(\alg{B}_{\geq m+1})$, then it
remains nilpotent in $H^{*}(\alg{B}_{\geq m})$.

By definition, $\alg{B}_{\geq m+1}$ has profile function
\[
\alg{B}_{\geq m+1}: (\underbrace{0, \dots, 0}_{m-1}, 0 , b_{m+1}, \dots).
\]
We define the algebra $\alg{B}_{\geq m+1,j}$ by the profile
function
\[
(\underbrace{0, \dots, 0}_{m-1}, j , b_{m+1}, \dots)
\]
for each $0 \leq j \leq b_{m}$, so $\alg{B}_{\geq m+1,0} =
\alg{B}_{\geq m+1} = \alg{B}_{\geq m+2, b_{m+1}}$. For each $j$ with
$0 \leq j \leq b_{m}-1$, we have a Hopf algebra extension
\[
\alg{B}_{\geq m+1,j} \to \alg{B}_{\geq m+1,j+1} \to D(\pst{j}{m}).
\]

\medskip
\begin{tikzpicture}[scale=0.5]
\draw (0,10) -- (0,0) -- (10,0);
\draw (4,0) -- (4,2) -- (5,2) -- (5,7) -- (6,7) -- (6,9);
\draw [dotted] (4,2) -- (4,3) -- (5,3);
\draw (4.5, 2.5) node {$\bullet$};
\draw (3,6) node {$\xi_{m}^{p^{j}}$};
\draw [->, shorten <=5mm, shorten >= 2mm] (3,6) -- (4.5, 2.5);
\draw[below] (4.5,0) node {\tiny{${m}$}};
\draw [right] (5,4.5) node {$\alg{B}_{\geq m+1,j}$};
\end{tikzpicture}

We proceed by induction, so we assume that $b_{nk}$ is nilpotent over
$\alg{B}_{\geq m+1,j}$, and we want to prove that $b_{nk}$ remains
nilpotent over $\alg{B}_{\geq m+1,j+1}$.

If $b_{nk}^{M}=0$ in the cohomology of $\alg{B}_{\geq m+1,j}$, then
$b_{nk}^{M}$ is in the kernel of the restriction map
$H^{*}(\alg{B}_{\geq m+1,j+1}) \to H^{*}(\alg{B}_{\geq m+1,j})$, so by
the Quillen-Venkov \autoref{lemma-quillen-venkov}, we see that
$b_{nk}^{2M} = b_{mj} x$ for some $x \in H^{*}(\alg{B}_{\geq
m+1,j+1})$. Now we appeal to the Nilpotence
\autoref{lemma-nilpotence}. We know that $m \leq n-1$ and $j < b_{m}
\leq \binom{m}{2}+1$, and given this, the lemma says that if $m+j+1
\leq k$ then $b_{mj} b_{nk}$ is nilpotent. We know that
$\binom{n}{2}+1 \leq k$, so we check that $m+j+1 \leq \binom{n}{2}+1$,
or equivalently that $m+j \leq \binom{n}{2}$. This is easy: since $m
\leq n-1$ and $j \leq \binom{m}{2}$, we have
\[
m+j \leq n-1 + \binom{n-1}{2} = \binom{n}{2}.
\]
Thus $b_{mj} b_{nk}$ is nilpotent. Since $b_{nk}^{2M}$ is a multiple
of $b_{mj}$, we conclude that $b_{nk}$ is nilpotent. This completes
the inductive step and hence the proof.
\end{proof}


\section{Restricting to the algebra \texorpdfstring{$\D$}{D}}\label{sec-D}

Fix a sub-Hopf algebra $\alg{C} \subseteq \A$.  The goal in this
section is to use \autoref{thm-reduction} to show that $\alg{C}
\cap \D$ detects nilpotence in $\stable{(\alg{C} \cap \Binom)}$ and
that $H^{*}(\alg{C} \cap \Binom) \to (H^{*}(\alg{C} \cap \D))^{\alg{C}
\cap \Binom}$ is an $F$-isomorphism. As in the previous section, we
work entirely with sub-Hopf algebras of $\P$ and then appeal to
\autoref{prop-P-to-A}.

Suppose that $\alg{B}$ is a sub-Hopf algebra of $\P$ with profile
function $(b_{1}, b_{2}, \dots)$. We say that $\alg{B}$ has a
\emph{finite profile function} if each $b_{i} < \infty$. Here
is the main theorem in this section:

\begin{theorem}\label{thm-bnn-finite-profile}
Suppose that $\alg{B} \subseteq \P$ has a finite profile function.
Given $k\geq n\geq 1$, if $\xi_{n}^{p^{k}}$ is primitive in
$\dual{B}$, then $b_{nk}$ is nilpotent in $H^{*}(\alg{B})$.
\end{theorem}

Given this, to apply \autoref{thm-reduction} and
\autoref{prop-P-to-A} to the inclusion $\D \subseteq \Binom$,
we intersect everything with $\P$ and then look at the Hopf algebra
with profile function which looks like $\D\cap \P$ for the first $n-1$
terms, and after that it looks like $\Binom \cap \P$: the profile
function has the form $(1, 2, 3, \dots, n-1, \binom{n}{2}+1,
\binom{n+1}{2}+1, \dots)$. We need to know whether the cohomology
element $b_{nn}$ is nilpotent over this Hopf algebra, but this follows
immediately from \autoref{thm-bnn-finite-profile}.

\begin{corollary}\label{cor-D}
For any $\alg{C} \subseteq \A$, $\alg{C} \cap \D$
detects nilpotence in $\cat{Stable}(\dual{(\alg{C} \cap \Binom)})$, and
the restriction map
\[
H^{*}(\alg{C} \cap \Binom) \to \left( H^{*}(\alg{C} \cap \D)
\right)^{\alg{C} \cap \Binom}
\]
is an $F$-isomorphism.
\end{corollary}

If we combine \cref{cor-binom,cor-D}, we get
\cref{thm-nilpotence,thm-f-iso}. The $F$-isomorphism
aspect is immediate: just compose the $F$-isomorphisms from the
corollaries. Detection of nilpotence perhaps require a bit more
consideration, so let's take care of that.

\begin{proof}[Proof of \autoref{thm-nilpotence}]
Given $\alg{C} \subseteq \A$, we know that $\alg{C} \cap \Binom$
detects nilpotence in $\stable{C}$, and also that $\alg{C} \cap \D$
detects nilpotence in $\cat{Stable}(\dual{(\alg{C} \cap \Binom)})$. We
may need to work in categories $\stable{B}$ for different sub-Hopf
algebras $\alg{B}$, so we may add a superscript to indicate the
algebra, for example writing $\pi_{*}^{\alg{B}}(-)$. Given a sub-Hopf
algebra $\alg{B} \subseteq \alg{C}$ and an object $X$ in $\stable{C}$,
we use the notation $X \restrict{\alg{B}}$ to denote the restriction
of $X$ to $\stable{B}$.

Detecting nilpotence is defined in terms of objects $H(\alg{C} \cap
\Binom)$ and $H(\alg{C} \cap \D)$, so for example we want to verify
this:
\begin{itemize}
\item Fix a ring spectrum $R$ in $\stable{C}$ and $\alpha \in
\pi_{*}^{\alg{C}}R$. Suppose that $H(\alg{C} \cap \D)_{*}\alpha =
0$. Then $\alpha$ is nilpotent.
\end{itemize}
Since $\alg{C} \cap \D$ detects nilpotence in
$\cat{Stable}(\dual{(\alg{C} \cap \Binom)})$, we know that $\alpha \in
\pi_{*}^{\alg{C} \cap \Binom}(R \restrict{\alg{C} \cap \Binom})$ is
nilpotent, and by replacing $\alpha$ with a power, we may assume that
this element is in fact zero. Now we appeal to \cite[Corollary
1.2.7]{palmieri;steenrod-stable}: 
\[
\pi_{*}^{\alg{C} \cap \Binom}(R \restrict{\alg{C} \cap \Binom}) \cong
\pi_{*}^{\alg{C}} (H(\alg{C} \cap \Binom) \Smash R) = 
H(\alg{C} \cap \Binom)_{*} R.
\]
Since $\alg{C} \cap \Binom$ detects nilpotence in $\stable{C}$, we
conclude that $\alpha \in \pi_{*}^{C} R$ is nilpotent.

One can similarly check the other conditions in the definition of
detecting nilpotence.
\end{proof}

To prove \autoref{thm-bnn-finite-profile}, we use the observation that
any $\alg{B} \subseteq \P$ with a finite profile function comes
equipped with surjections to finite-dimensional Hopf algebras: for
each $j \geq 1$ we have a quotient map $\alg{B} \to \alg{B}_{< j}$.
So if we can show that $\xi_{n}^{p^{k}}$ is primitive in some
$\alg{B}_{< j}$ and that $b_{nk} \in H^{*}(\alg{B}_{< j})$ is
nilpotent, then its image $b_{nk} \in H^{*}(\alg{B})$ is also
nilpotent. Since $\alg{B}_{< j}$ is finite-dimensional, we can check
nilpotence of $b_{nk}$ by appealing to \cite[Theorem
1.2]{palmieri;quasi}: we can detect nilpotence in $H^{*}(\alg{B}_{<
j})$ by restricting to the ``quasi-elementary'' sub-Hopf algebras of
$\alg{B}_{< j}$. This means that we need to understand, at least
partially, the quasi-elementary sub-Hopf algebras of an algebra like
$\alg{B}_{< j}$ --- that is, a finite-dimensional sub-Hopf algebra of
$\P_{< j}$. We investigate this in the next subsection.


\subsection{Quasi-elementary Hopf algebras}

Serre showed that a $p$-group is elementary abelian if and only if its
mod $p$ cohomology satisfies a certain condition, and that
cohomological condition is a key part in proving many facts about
group cohomology, such as detecting nilpotence and Quillen
stratification.  Wilkerson coined the term ``quasi-elementary'' for
those Hopf algebras which satisfy Serre's cohomological criterion
(Definition~\ref{defn-quasi}). It is clear that any elementary Hopf
algebra is quasi-elementary, and he showed that the converse holds for
sub-Hopf algebras of the mod 2 Steenrod algebra: this is very useful
in proving results at the prime 2 like
\autoref{thm-2-elementary}. Unfortunately these concepts differ
for sub-Hopf algebras of the odd primary Steenrod algebra: there are
quasi-elementary sub-Hopf algebras of $\A$ which are not
elementary. (The same goes for sub-Hopf algebras of $\P$: the issue is
not the exterior part of the Steenrod algebra.)  The main obstacle to
obtaining a result like \autoref{thm-2-elementary} at odd primes
is that we do not have a full classification of the quasi-elementary
sub-Hopf algebras of $\A$ or of $\P$. Nakano and the author claimed in
\cite{nakano-palmieri} to have such a classification. The proof is
incomplete, though, and we do not know how to fill in the gaps. We
provide a partial replacement here: see \autoref{thm-quasi}
below. We also discuss the classification problem in
\autoref{sec-discussion}.

Serre's cohomological condition says that a $p$-group is elementary
abelian if and only if every product of Bocksteins applied to nonzero
cohomology classes in $H^{1}$ is non-nilpotent. For Hopf algebras, it
is often not necessary to look at all degree one elements --- this is
true for finite-dimensional Hopf algebras, and for sub-Hopf algebras
of $\A$, it is also true for those with finite profile function.
Indeed if we also suppose that we are in the evenly-graded connected
case, then it suffices to look at degree one elements in the kernel of
the first nonvanishing Steenrod operation on $H^{s,t}$, which is
$\PP{t/2}$ or $\tildePP{0}$, depending on the grading
convention. Since we are only dealing with evenly graded Hopf algebras
in this section, we will use the $\widetilde{\mathcal{P}}$ notation.

The following definition was first published in \cite{palmieri;quasi},
although as noted above, the idea is really due to Wilkerson. We focus
on the case of evenly graded connected Hopf algebras, and that
simplifies the definition a bit.

\begin{definition}\label{defn-quasi}
Fix an odd prime $p$, and let $\alg{Q}$ be a connected evenly graded
cocommutative Hopf algebra over $\Fp$.  A nonzero element $v \in
H^{2}(\alg{Q})$ is a \emph{Serre element} if $v=\beta \tildePP{0} (w)$
for some $w \in H^{1}(\alg{Q}) \cap \ker \tildePP{0}$.  The Hopf
algebra $\alg{Q}$ is {\em quasi-elementary} if no product of Serre
elements is nilpotent.
\end{definition}

Note first that every elementary Hopf algebra is quasi-elementary.
Note also that, as perhaps suggested above, the definition of
quasi-elementary does not work well for sub-Hopf algebras of $\P$ with
infinite profile functions, or indeed for $\P$ itself: there are no
Serre elements in $H^{*}(\P)$, for example, because $\tildePP{0}$ is
injective on $H^{1}(\P)$. So in the case of sub-Hopf algebras of $\P$,
we should add conditions, for example assuming a finite profile
function, in order to get meaningful results.

Here is a partial classification of quasi-elementary sub-Hopf algebras
of $\P$; this is the main result of this section.

\begin{theorem}\label{thm-quasi}
Fix $N$ with $1 < N \leq \infty$ and suppose that $\alg{Q} \subseteq
\P_{< N}$ has finite profile function $(q_{1}, q_{2}, q_{3},
\dots)$. If $\alg{Q}$ is quasi-elementary, then $q_{i} \leq i$ when $i
< N/2$. If $N = \infty$, then $\alg{Q}$ has a bounded profile
function.
\end{theorem}

(Recall from \autoref{remark-adams-margolis-A<n} that the sub-Hopf
algebras of $\P_{< N}$ are determined by their profile functions, just
as for sub-Hopf algebras of $\P$.)

\begin{corollary}
Every quasi-elementary sub-Hopf algebra of $\P$ is contained in $\D$
and has bounded profile function.
\end{corollary}

This corollary verifies \cite[Conjecture
2.1.13]{palmieri;steenrod-stable}. See \autoref{sec-discussion}
for discussion regarding the quasi-elementary sub-Hopf algebras of
$\P$.

The following is a useful criterion for identifying
non-quasi-elementary sub-Hopf algebras of $\P_{< N}$.

\begin{lemma}\label{lemma-quasi}
Fix $N$ with $1 < N \leq \infty$. Fix integers $t,n <N$, and fix a
sub-Hopf algebra $\alg{Q}$ of $\P_{< N}$ with profile function
$(q_{1}, q_{2}, \dots)$. Suppose that
\begin{enumerate}
\item $\xi_{t}^{p^{s}}$ and $\xi_{n}^{p^{k}}$ are both nonzero
and primitive in $\dual{Q}$, and
\item The profile function entries $q_{t}$ and $q_{n}$ are finite.
\end{enumerate}
If $b_{ts} b_{nk}$ is nilpotent in $H^{*}(\alg{Q})$, then $\alg{Q}$ is
not quasi-elementary.
\end{lemma}

Of course condition (b) is automatically satisfied if $\alg{Q}$ has a
finite profile function. Note that the integers $t,n$ need not be
distinct.

\begin{proof}
Let $L = q_{t}-1$ and $M=q_{n}-1$, so that $\xi_{t}^{p^{L}}$ and
$\xi_{n}^{p^{M}}$ are the largest $p$th powers of these generators
which are nonzero in $\dual{Q}$. Then $b_{t,L}$ and $b_{n,M}$ are
Serre elements. (Since $\xi_{t}^{p^{s}}$ and $\xi_{n}^{p^{k}}$ are
primitive, so are any larger $p$th powers of these generators, and so
$b_{t,L}$ and $b_{n,M}$ are elements of $H^{2}$.)  We will show that
$b_{t,L} b_{n,M}$ is nilpotent.

Suppose that $(b_{ts} b_{nk})^{i}=0$. We may apply the algebra
map $\tildePP{0}$ repeatedly to get $b_{t,s+\ell}^{i} b_{n,k+\ell}^{i}=0$
for any $\ell \geq 0$, and so without loss of generality, we may
assume that $k=M$: we may assume that the relation has the form
$b_{t,s}^{i} b_{nM}^{i} = 0$, with $s \leq L$. If $s = L$, then we
have a product $b_{tL}b_{nM}$ of Serre elements which is nilpotent, so
we may assume that $s<L$.

Choose $d$ so that $p^{d}>i$, and multiply this relation by
$b_{nM}^{p^{d}-i}$ to get $b_{ts}^{i}b_{nM}^{p^{d}} = 0$. Now apply
Steenrod operations that increase the power of $b_{nM}$ and increase
the second index of $b_{ts}$. That is,
apply $\tildePP{p^{d+j-1}} \cdots \tildePP{p^{d+1}}
\tildePP{p^{d}}$ to get
\[
b_{t,s+j}^{i}b_{nM}^{p^{d+j}} = 0.
\]
When $j=L-s$, we get $b_{t,L}^{i} b_{nM}^{p^{d+L-s}} = 0$, as desired.
\end{proof}

Now we use the Nilpotence \autoref{lemma-nilpotence} to prove
\autoref{thm-quasi}: if $\alg{Q}$ does not satisfy the hypotheses,
we find relations $b_{t,s}^{N} b_{n,k}^{N}$ in $H^{*}(\alg{Q})$ and
deduce from \autoref{lemma-quasi} that $\alg{Q}$ is not
quasi-elementary.

\begin{proof}[Proof of \autoref{thm-quasi}]
Fix a sub-Hopf algebra $\alg{Q}$ of $\P_{<N}$ with finite profile
function.

If $N=\infty$ and if the profile function for $\alg{Q}$ is not bounded,
then find $t$ so that $\xi_{i}=0$ for $i<t$ while $\xi_{t} \neq
0$. Let $s$ be the largest number so that $\xi_{t}^{p^{s}} \neq 0$.
Since the profile function is not bounded, we can find $k \geq s+t+1$
so that $\xi_{n}^{p^{k}} \neq 0$ for some $n$. Find the minimal such
$n$ and apply \cref{lemma-nilpotence,lemma-quasi} to
deduce that $\alg{Q}$ is not quasi-elementary. This takes care of the
last sentence of \autoref{thm-quasi}.

Now assume that $N$ is arbitrary, and assume that some entry in the
profile function satisfies $q_{n} > n$ with $n<N/2$. Choose the
minimal such $n$, meaning that $\xi_{j}^{p^{j}}=0$ for all $j<n$. We
want to show that $\alg{Q}$ is not quasi-elementary. If $\xi_{n}$ is
the first non-vanishing generator --- that is, if $\xi_{i}=0$ for all
$i<n$ --- then essentially by the argument in the proof of
\cite[Proposition 4.1]{miller-wilkerson}, we can see that $\alg{Q}$ is
not quasi-elementary.

In more detail (to fill in what is meant by ``essentially'' in the
previous sentence): if $\xi_{i}=0$ for all $i<n$ and if
$\xi_{n}^{p^{n}} \neq 0$, then the reduced coproduct on $\xi_{2n}$ is
$\xi_{n}^{p^{n}} \otimes \xi_{n}$. Since this tensor product is
nonzero in $\dual{Q} \otimes \dual{Q}$, $\xi_{2n}$ must be nonzero in
$\dual{Q}$; equivalently, $\xi_{2n} \neq 0$ can be deduced from
$\xi_{n}^{p^{n}} \neq 0$ and the Adams-Margolis
\autoref{thm-adams-margolis} on profile functions.  The reduced coproduct
$\xi_{2n} \mapsto \xi_{n}^{p^{n}} \otimes \xi_{n}$ yields the relation
$h_{nn} h_{n0} = 0$ in cohomology (both $h_{n0}$ and $h_{nn}$ are
cohomology classes because $\xi_{n}$ is primitive in
$\dual{Q}$). Applying Steenrod operations to this relation gives
\begin{align*}
\beta \tildePP{0}: & \quad b_{nn}h_{n1}  - h_{n,n+1} b_{n0} = 0, \\
\beta \tildePP{1}: & \quad b_{nn}^{p} b_{n1}  - b_{n,n+1} b_{n0}^{p} = 0, \\
\beta \tildePP{p}: & \quad b_{nn}^{p^{2}} b_{n2} - b_{n,n+2} b_{n0}^{p^{2}} = 0,
\end{align*}
and in general 
\[
b_{nn}^{p^{d}} b_{nd} - b_{n,n+d} b_{n0}^{p^{d}} = 0.
\]
By assumption $\alg{Q}$ has a finite profile function, so
$\xi_{n}^{p^{n+d}}=0$ for some $d$, and if we choose the smallest $d$
making this hold, then we get the monomial relation $b_{nn}^{p^{d}}
b_{nd}=0$. By \autoref{lemma-quasi}, $\alg{Q}$ is not
quasi-elementary. This deals with the case when $\xi_{n}^{p^{n}} \neq
0$ and $\xi_{i}=0$ for $i<n$.

To deal with the remaining case, we may assume that for some $t<n$, we
have in $\dual{Q}$:
\begin{enumerate}
\item[(1)] $\xi_{i} = 0$ for $i<t$,
\item[(2)] $\xi_{t} \neq 0$,
\item[(3)] $\xi_{j}^{p^{j}}=0$ for all $j < n$,
\item[(4)] $\xi_{n}^{p^{n}} \neq 0$.
\item[(5)] $\xi_{k}^{p^{j}}=0$ for all $j,k$ with $N> k+t \geq 2t+1$ and $j \geq 2t$.
\end{enumerate}
Explanation: (1) and (2) say that $\xi_{t}$ is the first generator
present in $\dual{Q}$. (3) and (4) say that $\xi_{n}$ is the first
generator with $\xi_{n}^{p^{n}} \neq 0$.  (5)
is because of annihilator considerations: if $\xi_{t}^{p^{s}} \neq 0$
with $s<t$ and $\xi_{k}^{p^{j}} \neq 0$, then by
\autoref{lemma-nilpotence}, $b_{ts} b_{kj}$ is nilpotent if
$N>k+t$, $k \geq t+1$, and $j \geq 2t \geq t+s+1$, so if any such
$\xi_{k}^{p^{j}}$ were nonzero, we would get a monomial relation in
$H^{*}$, and \autoref{lemma-quasi} would then tell us that $\alg{Q}$
is not quasi-elementary.

Combining (4) and (5), along with the assumption that $n>t$, gives
\begin{enumerate}
\item [(6)] $n < 2t$, or equivalently $n \leq 2t-1$.
\end{enumerate}
That is: $n$ satisfies $n>t$ so $n \geq t+1$, so $n+t \geq 2t+1$. Also
$2n<N$ so $n+t<2n<N$, so $n$ satisfies the first inequality in (5): $N
> n+t \geq 2t+1$. Since $\xi_{n}^{p^{n}}\neq 0$, the exponent $n$ must fail to
satisfy the second inequality in (5): we must have $n < 2t$.

The reduced coproduct on $\xi_{2n}$ is
\[
\sum_{j=1}^{2n-1} \xi_{2n-j}^{p^{j}} \otimes \xi_{j}.
\]
Because of (1), we can change the limits on the sum to go from $t$ to
$2n-t$. Because of (3), we can omit more terms: when $j>n$, then
$2n-j<n<j$, so $\xi_{2n-j}^{p^{j}}=0$. So the reduced coproduct on
$\xi_{2n}$ is
\[
\sum_{j=t}^{n} \xi_{2n-j}^{p^{j}} \otimes \xi_{j}.
\]

\begin{sublemma}
If the Hopf algebra $\dual{Q}$ satisfies conditions (1)--(6), then
the elements $\xi_{2n-j}^{p^{j}}$ and $\xi_{j}$ are primitive when $t
\leq j \leq n$.
\end{sublemma}

\begin{proof}
Combining (6) with the inequality $j \leq n$ gives $j< 2t$. Therefore
each element $\xi_{j}$ is primitive: each term in its coproduct
involves $\xi_{i}$ and $\xi_{j-i}$, and either $i$ or $j-i$ will be
less than $t$, so assumption (1) tells us that each term is zero.

Now we examine $\xi_{2n-j}^{p^{j}}$ when $t \leq j \leq n$. The
reduced coproduct on this element is
\[
\xi_{2n-j}^{p^{j}} \mapsto \sum_{i=1}^{2n-j-1} \xi_{2n-j-i}^{p^{i+j}}
\otimes \xi_{i}^{p^{j}}.
\]
The terms with $i<t$ are zero, by (1). Combining (6) with $i \geq t$
and $j \geq t$, we get
\[
n \geq n + (n-2t + 1) = 2n - 2t + 1 \geq 2n - j - i + 1 > 2n - j - i.
\]
This means that assumption (3) applies to $\xi_{2n-j-i}$.

We also have
\[
2i+2j \geq 2t + 2t = 4t > 2n,
\]
so
\[
i+j > 2n - j - i.
\]
So by (3), we have $\xi_{2n-j-i}^{p^{i+j}} = 0$ for all terms in the
coproduct. Therefore $\xi_{2n-j}^{p^{j}}$ is primitive.
\end{proof}

As a result, the coproduct on $\xi_{2n}$ produces a relation
\[
\sum_{j=t}^{n} h_{2n-j,j} h_{j0} = 0
\]
where each individual factor $h_{k\ell}$ is a cohomology class, so we can
apply Steenrod operations and use the Cartan formula. Applying $(\beta
\tildePP{1}) (\beta \tildePP{0})$ to this yields
\[
\sum_{j=t}^{n} (b_{2n-j,j}^{p} b_{j1} - b_{2n-j,j+1} b_{j0}^{p}) = 0.
\]
Now apply $\tildePP{p^{n-2}} \cdots \tildePP{p^{2}} \tildePP{p}$:
\[
\sum_{j=t}^{n} (b_{2n-j,j}^{p^{n-1}} b_{j,n-1} - b_{2n-j,n+j-1} b_{j0}^{p^{n-1}}) = 0.
\]
If $j<n$, then $b_{j,n-1}=0$ by (3). Also $2t \leq n+t-1 \leq n+j-1$,
so $b_{2n-j,n+j-1}=0$ by (5). So all terms with $j<n$ vanish.

That leaves us with the two $j=n$ terms, the second of which involves
$b_{n,2n-1}$, but $2n-1 \geq 2t$, so $b_{n,2n-1}=0$ by (5). So we have a
monomial relation: $b_{n,n}^{p^{n-1}} b_{n,n-1} = 0$, and therefore
$\alg{Q}$ cannot be quasi-elementary by \autoref{lemma-quasi}. Note
that applying $\tildePP{p^{n-1}}$ to this yields $b_{n,n}^{p^{n}+1} = 0$,
if you want a ``cleaner'' relation.
\end{proof}


\subsection{Proof of \autoref{thm-bnn-finite-profile}}

Now we prove \autoref{thm-bnn-finite-profile}, so suppose that
$\alg{B} \subseteq \P$ has a finite profile function. Assume that $k
\geq n$ and that $\xi_{n}^{p^{k}}$ is primitive and nonzero in
$\dual{B}$. We want to show that $b_{nk}$ is nilpotent in
$H^{*}(\alg{B})$.

There is a quotient map $\alg{B} \to \alg{B}_{< j}$ for each $j >0$,
and since $\alg{B}$ has a finite profile function, $\alg{B}_{<j}$ is
finite-dimensional.  Note that if $j>n$, then $b_{nk}$ is a cohomology
class in $H^{*}(\alg{B}_{<j})$, and $b_{nk} \in H^{*}(\alg{B}_{< j})$
maps to $b_{nk} \in H^{*}(\alg{B})$.  If we can show that $b_{nk}$ is
nilpotent in the cohomology of $H^{*}(\alg{B}_{< j})$ for some $j>n$,
the proof will be complete. 

Choose finite $N$ so that $N > 2n$. By replacing $\alg{B}$ by
$\alg{B}_{<N}$, we may assume that we have a finite-dimensional
sub-Hopf algebra $\alg{B} \subseteq \P_{< N}$ such that
$\xi_{n}^{p^{k}}$ is primitive in $\dual{B}$ and that $k\geq n$. Now
we are in position to use \autoref{thm-quasi}, the partial
classification of quasi-elementary sub-Hopf algebras in $\P_{< N}$: in
the dual of every quasi-elementary sub-Hopf algebra $\alg{Q} $ of
$\alg{B}$, we have $\xi_{n}^{p^{n}} = 0$, so (since $k\geq n$)
$\xi_{n}^{p^{k}}=0$, so $b_{nk} \mapsto 0$ via the restriction map
$H^{*}(\alg{B}) \to H^{*}(\alg{Q})$. Since quasi-elementary Hopf
algebras detect nilpotence in finite-dimensional Hopf algebras by
\cite[Theorem 1.2]{palmieri;quasi}, we conclude that $b_{nk} \in
H^{*}(\alg{B})$ is nilpotent.  \qed


\section{Bounded profile functions}\label{sec-bounded}

At the prime 2, the elementary sub-Hopf algebras of the mod 2 Steenrod
algebra $\A$ detect nilpotence and give an $F$-isomorphism $H^{*}(\A)
\to \left( \varprojlim H^{*}(\alg{E}) \right)^{\A}$ --- see
\autoref{thm-2-elementary}.  The goal of this section is to prove
\autoref{thm-bdd-nilpotence}, a partial replacement for this.

Recall the theorem: given a sub-Hopf algebra $\alg{C} \subseteq \A$,
let $\alg{C}' = \alg{C} \cap \D$ and let $\Bdd(\alg{C}')$ be the
category whose objects are sub-Hopf algebras $\alg{B} \subseteq
\alg{C}'$ with bounded profile function and whose morphisms are
inclusions. The theorem asserts that the objects in $\Bdd(\alg{C}')$
detect nilpotence in $\stable{C}$. We will get to the proof of this
soon, but first let's investigate the potential of having an
$F$-isomorphism in this setting.  We note that we can define an action
of $\alg{C}$ on $\varprojlim_{\Bdd(\alg{C}')} H^{*}(\alg{B})$ as
follows: let $\nBdd(\alg{C}')$ be the full subcategory of
$\Bdd(\alg{C}')$ consisting of the normal sub-Hopf algebras of
$\alg{C}'$ with bounded profile function. This is cofinal in
$\Bdd(\alg{C}')$, so we have
\[
\varprojlim_{\Bdd(\alg{C}')} H^{*}(\alg{B}) \cong
\varprojlim_{\nBdd(\alg{C}')} H^{*}(\alg{B}).
\]
There are compatible actions of $\alg{C}$ on each term on the right,
and therefore an induced action on the limit. For $\alg{B} \in \nBdd
(\alg{C}')$, the map $H^{*}(\alg{C}') \to H^{*}(\alg{B})$ is one of
$\alg{C}$-algebras, and so we have an induced map
\[
\left( H^{*}(\alg{C}') \right)^{\alg{C}} \to
\left(\varprojlim_{\Bdd(\alg{C}')} H^{*}(\alg{B})\right)^{\alg{C}}.
\]
\autoref{thm-f-iso} says that we have an $F$-isomorphism
$H^{*}(\alg{C}) \to \left( H^{*}(\alg{C}') \right)^{\alg{C}}$.  Once
we have established that $\Bdd(\alg{C}')$ detects nilpotence in
$\stable{C}$, to show that the map
\[
H^{*}(\alg{C}')^{\alg{C}} \to \left( \varprojlim_{\Bdd (\alg{C}')}
H^{*}(\alg{B}) \right)^{\alg{C}}
\]
is an $F$-isomorphism, we would just need to show that this is an
$F$-epimorphism. In fact it suffices to show that the map
\[
H^{*}(\alg{C}') \to \varprojlim_{\Bdd (\alg{C}')} H^{*}(\alg{B})
\]
is an $F$-epimorphism.  At the prime 2 the corresponding result, in
which the inverse limit is taken over the elementary sub-Hopf algebras
of $\alg{C}$, can be established by combining several facts: first,
one can compute the cohomology of the elementary sub-Hopf algebras and
the induced maps on cohomology; second, because every Hopf algebra
involved here has a finite profile function, it is an inverse limit of
finite-dimensional Hopf algebras; and third, a theorem of Hopkins and
Smith \cite[Theorem 4.13]{hopkins-smith} says that if $\alg{\Lambda}
\subseteq \alg{\Gamma}$ is an inclusion of finite-dimensional graded
connected Hopf algebras, then the induced map on cohomology is an
$F$-epimorphism. See \cite[Theorem 5.3.4]{palmieri;steenrod-stable}
for details.

In our situation we have all of these components except for the
computation of the cohomology of each algebra. If we could do that, we
presumably would get an $F$-isomorphism. The cohomology of an
arbitrary object in $\Bdd (\alg{C}')$ is likely to be quite
complicated, but we hope that it would be good enough to work with the
quasi-elementary sub-Hopf algebras of $\alg{C}$. Recall from
\autoref{thm-quasi} that every quasi-elementary sub-Hopf algebra of
$\A$ is contained in $\D$ and has a bounded profile function, so the
category $\Bdd$ contains the quasi-elementary sub-Hopf
algebras. (Note, though, that there are Hopf algebras with 
bounded profile function that are not quasi-elementary, such as the
Hopf algebra with profile function $(1,2,3,3,3,\dots)$.)

If \autoref{conj-quasi} holds, then we would get an $F$-epimorphism
\[
H^{*}(\alg{C}') \to \varprojlim_{\cat{QElem}(\alg{C}')} H^{*}(\alg{Q}),
\]
where the inverse limit is taken over the quasi-elementary sub-Hopf
algebras of $\alg{C}'$. Given the conjecture, we would also likely get
a version of \autoref{thm-bdd-nilpotence}: see the last paragraph of
this section. So we are led to this conjecture:

\begin{conjecture}\label{conj-f-iso-quasi}
For $\alg{C} \subseteq \A$, the quasi-elementary sub-Hopf
algebras of $\alg{C}$ detect nilpotence in $\stable{C}$, and
there is an $F$-isomorphism
\[
H^{*}(\alg{C}) \to \left( \varprojlim H^{*}(\alg{Q}) \right)^{\alg{C}},
\]
where the limit is taken over the category whose objects are
quasi-elementary sub-Hopf algebras $\alg{Q}$ of $\alg{C}$, and whose
morphisms are inclusions.
\end{conjecture}

Anyway, back to our task, proving \autoref{thm-bdd-nilpotence}. We
do this by imitating the $p=2$ proof from \cite[Section
5.3]{palmieri;steenrod-stable}. We work in the category
$\stable{C}$. Let's start by recalling
\cite[1.5.1--2]{palmieri;steenrod-stable}.

\begin{lemma}\label{lemma-bousfield-ravenel}
\begin{enumerate}
\item Given a map $f: X \to X$ (omitting suspensions) with
cofiber $X/f$ and telescope $f^{-1}X$, we have $\bousfield{X} =
\bousfield{X/f} \vee \bousfield{f^{-1} X}$ and $X/f \Smash f^{-1} X = 0$.
\item Therefore given Hopf algebras $\alg{B}' \subseteq \alg{B}
\subseteq \alg{C}$ with $\alg{B}'$ normal in $\alg{B}$ and $\alg{B} \modmod
\alg{B}' \cong D(\pst{s}{t})$, we have $\bousfield{\H{B}} =
\bousfield{\H{B}'} \vee \bousfield{b_{st}^{-1} \H{B}}$ and
$\H{B}' \Smash b_{st}^{-1} \H{B} = 0$.
\item Similarly, given Hopf algebras $\alg{B}' \subseteq \alg{B}
\subseteq \alg{C}$ with $\alg{B}'$ normal in $\alg{B}$ and $\alg{B} \modmod
\alg{B}' \cong \exterior (Q_{n})$, we have $\bousfield{\H{B}} =
\bousfield{\H{B}'} \vee \bousfield{a_{n}^{-1} \H{B}}$ and
$\H{B}' \Smash a_{n}^{-1} \H{B} = 0$.
\end{enumerate}
\end{lemma}

In the ensuing discussion, we will often start with a Hopf algebra
quotient map $\alg{B} \to D(\pst{s}{t})$; note that this is equivalent
to having a Hopf algebra extension $\alg{B}' \to \alg{B} \to
D(\pst{s}{t})$, and similarly for a quotient map $\alg{B} \to
\exterior (Q_{n})$. We abuse notation and write $\alg{B}/b_{ts}$ for
the Hopf algebra kernel: this is the sub-Hopf algebra of $\alg{B}$
whose dual is $\dual{B}/(\xi_{t}^{p^{s}})$. According to
\cite[1.5.2]{palmieri;steenrod-stable}, $H(\alg{B}/b_{ts})$ is not the
cofiber of $\H{B} \xrightarrow{b_{ts}} \H{B}$, that is,
$H(\alg{B}/b_{ts}) \neq (\H{B})/b_{ts}$. They are however Bousfield
equivalent, and that's all we will need below. As a result, in
Bousfield class computations we will omit parentheses, writing
$\bousfield{H\alg{B}/b_{ts}}$ for $\bousfield{H(\alg{B}/b_{ts})} =
\bousfield{(H\alg{B})/b_{ts}}$.  Similarly we write $\alg{B}/a_{n}$
for the algebra dual to $\dual{B}/(\tau_{n})$, and for this one we do
have $H(\alg{B}/a_{n}) = (\H{B})/a_{n}$.

For example, $\alg{B}'$ in part (b) above could be written as
$\alg{B}/b_{ts}$, and we have
\[
\bousfield{\H{B}} = \bousfield{\H{B}/b_{ts}} \vee
\bousfield{b_{ts}^{-1} \H{B}}.
\]

We also define $\alg{B}/0 := \alg{B}$, so $\alg{B}/(0, x, y, z,
\dots) = \alg{B}/(x, y, z, \dots)$. This convention will be useful in
\autoref{lemma-bousfield-big}.

\begin{corollary}\label{cor-bousfield}
Suppose that $\alg{B} \subseteq \alg{C}$ has profile function
$((b_{i});(e_{j}))$
\begin{enumerate}
\item Suppose that $b_{r}, b_{q} < \infty$ for some $r, q$, and let
$i=b_{r}-1$, $j=b_{q}-1$. Assume that we have Hopf algebra quotient
maps $\alg{B} \to D(\pst{i}{r})$ and $\alg{B} \to D(\pst{j}{q})$. Then
we have
\[
\bousfield{\H{B}} = \bousfield{\H{B} / b_{ri}}
\vee \bousfield{b_{ri}^{-1} \H{B}/b_{qj}}
\vee \bousfield{b_{qj}^{-1} b_{ri}^{-1} \H{B}}.
\]
\item Suppose that $b_{r} < \infty$ and $e_{n} = 2$ for some $r, n$,
and let $i=b_{r}-1$. Assume that we have Hopf algebra quotient maps
$\alg{B} \to D(\pst{i}{r})$ and $\alg{B} \to \exterior (Q_{n})$. Then
we have
\[
\bousfield{\H{B}} = \bousfield{\H{B} / b_{ri}}
\vee \bousfield{b_{ri}^{-1} \H{B}/a_{n}} 
\vee \bousfield{a_{n}^{-1} b_{ri}^{-1} \H{B}}.
\]
\end{enumerate}
\end{corollary}

Of course if $b_{qj} b_{ri}$ is nilpotent, the last term in part (a)
drops out, and similarly if $a_{n} b_{ri}$ is nilpotent for part (b). 
This is how we will apply the result.

\begin{proof}
The lemma tells us that
\[
\bousfield{\H{B}} = \bousfield{\H{B} / b_{ri}}
\vee \bousfield{b_{ri}^{-1} \H{B}}.
\]
For part (a), the lemma also tells us that 
\[
\bousfield{\H{B}} = \bousfield{\H{B} / b_{qj}}
\vee \bousfield{b_{qj}^{-1} \H{B}}.
\]
Smash this with $\bousfield{b_{r,i}^{-1} \H{B}}$ and use
$\bousfield{\H{B} \Smash \H{B}} = \bousfield{\H{B}}$ to get
\[
\bousfield{b_{ri}^{-1} \H{B}} = \bousfield{b_{ri}^{-1} \H{B} / b_{qj}}
\vee \bousfield{b_{ri}^{-1} b_{q,j}^{-1} \H{B}}.
\]
Part (b) is similar.
\end{proof}

Here are the main players in the proof of
\autoref{thm-bdd-nilpotence}. Recall from \autoref{subsec-hopf}
that for a non-negative integer $r$, $\D_{\geq r}$ has profile function
\[
((\underbrace{0, \dots, 0}_{r-1}, r, r+1, r+2,\dots); 
(\underbrace{1, \dots, 1}_{r}, 2, 2, 2, \dots)).
\]
If $i \geq 0$ and $q \geq r+i$, we let $\D_{\geq r,i,q}$
be the sub-Hopf algebra of $\A$ with profile function
\[
((\underbrace{0, \dots, 0}_{r-1}, 
\underbrace{r, r+1, \dots, r+i}_{i+1},
\underbrace{r+i+1, \dots, r+i+1}_{q-(r+i+1)}, q, q+1, q+2, \dots);
(\underbrace{1, \dots, 1}_{r}, 2, 2, 2, \dots)).
\]
If $i \geq 0$, we define $\alg{B}_{r,i}$ by the (bounded) profile
function
\[
((\underbrace{0, \dots, 0}_{r-1}, 
\underbrace{r, r+1, \dots, r+i}_{i+1}, r+i+1,
r+i+1, r+i+1, \dots);
(\underbrace{1, \dots, 1}_{r}, 2, 2, 2, \dots)).
\]
Here are pictures of the polynomial components of the profile
functions for $\D_{\geq r,i,q}$ and $\alg{B}_{r,i}$:

\medskip
\begin{tikzpicture}[scale=0.5]
\draw (0,10) -- (0,0) -- (10,0);
\draw (2,0) -- (2,2) -- (4,4) -- (7,4) -- (7,7) -- (10,10);
\draw [dotted] (4,0) -- (4,4) -- (0,4);
\draw [dotted] (7,0) -- (7,7);
\draw[below] (2.5,0) node {\tiny{$_{r}$}};
\draw[below] (4.5,0) node {\tiny{$_{r+i+1}$}};
\draw[below] (7.5,0) node {\tiny{$_{q}$}};
\draw[left] (0,4.2) node {\tiny$_{r+i+1}$};
\draw[above] (5, 4) node {$\D_{\geq r,i,q}$};
\end{tikzpicture}
\hspace{0.2in}
\begin{tikzpicture}[scale=0.5]
\draw (0,10) -- (0,0) -- (10,0);
\draw (2,0) -- (2,2) -- (4,4) -- (10,4);
\draw [dotted] (4,0) -- (4,4) -- (0,4);
\draw[below] (2.5,0) node {\tiny{$_{r}$}};
\draw[below] (4.5,0) node {\tiny{$_{r+i+1}$}};
\draw[left] (0,4.2) node {\tiny$_{r+i+1}$};
\draw[above] (7, 4) node {$\alg{B}_{r,i}$};
\end{tikzpicture}

\begin{lemma}\label{lemma-colimit-hdr}
\begin{enumerate}
\item In $\stable{A}$ we have $\varinjlim_{r} \H{\D}_{\geq r} = H\Fp$,
and for fixed $r,i$, $\varinjlim_{q} \H{\D}_{\geq r,i,q} = \H{B}_{r,i}$. 
\item Given $\alg{C} \subseteq \A$, then in $\stable{C}$ we have
$\varinjlim_{r} H(\alg{C} \cap \D_{\geq r}) = H\Fp$, and for fixed $r,i$,
$\varinjlim_{q} H(\alg{C} \cap \D_{\geq r,i,q}) = H(\alg{C} \cap \alg{B}_{r,i})$.
\end{enumerate}
\end{lemma}

\begin{proof}
This is immediate. (See \cite[Lemma 5.3.2]{palmieri;steenrod-stable}
for the $p=2$ version, with a similarly detailed proof.)
\end{proof}

\begin{lemma}\label{lemma-bousfield-hdr}
Fix $r \geq 0$.
\begin{enumerate}
\item In $\stable{A}$ we have
\[
\bousfield{\H{D}_{\geq r}} = \bousfield{a_{r}^{-1} \H{D}_{\geq r}} \vee 
\bigvee_{i=0}^{r-1} \bousfield{b_{ri}^{-1} \H{D}_{\geq r} / (a_{r},
b_{r0}, \dots, b_{r,i-1})} \vee 
\bousfield{\H{D}_{\geq r+1}}.
\]
\item Given $\alg{C} \subseteq \A$ with profile function $((c_{i});
(e_{j}))$, in $\stable{C}$ we have a similar formula for
$\bousfield{H(\alg{C} \cap \D_{\geq r})}$, but it depends on whether
$\tau_{r}$ is zero in $\dual{C}$:
\begin{itemize}
\item If $\tau_{r} = 0$ in $\dual{C}$, then
\[
\bousfield{H(\alg{C} \cap \D_{\geq r})} = 
\bigvee_{i=0}^{\min (c_{r}, r)-1} \bousfield{b_{ri}^{-1} H(\alg{C}
\cap \D_{\geq r}) / (b_{r0}, \dots, b_{r,i-1})} \vee 
\bousfield{H(\alg{C} \cap \D_{\geq r+1})}.
\]
\item If $\tau_{r} \neq 0$ in $\dual{C}$, then
\begin{align*}
\bousfield{H(\alg{C} \cap \D_{\geq r})} = &
\bousfield{a_{r}^{-1} H(\alg{C} \cap \D_{\geq r})} \\
&\vee 
\bigvee_{i=0}^{\min (c_{r}, r)-1} \bousfield{b_{ri}^{-1} H(\alg{C} \cap \D_{\geq r}) / (a_{r},
b_{r0}, \dots, b_{r,i-1})} \\
& \vee \bousfield{H(\alg{C} \cap \D_{\geq r+1})}
\end{align*}
\end{itemize}
\end{enumerate}
\end{lemma}

\begin{proof}
(a) By the previous lemma we have
\[
\bousfield{\H{D}_{\geq r}} = 
\bousfield{a_{r}^{-1} \H{D}_{\geq r}} \vee \bousfield{\H{D}_{\geq r}/a_{r}}.
\]
Now we apply the lemma repeatedly to the last term: we find that
\begin{align*}
\bousfield{\H{D}_{\geq r}/a_{r}} &= 
\bousfield{b_{r0}^{-1} \H{D}_{\geq r}/a_{r}} \vee 
\bousfield{\H{D}_{\geq r}/(a_{r}, b_{r0})} \\
&=
\bousfield{b_{r0}^{-1} \H{D}_{\geq r}/a_{r}} \vee 
\bousfield{b_{r1}^{-1} \H{D}_{\geq r}/(a_{r}, b_{r0})} \vee 
\bousfield{\H{D}_{\geq r}/(a_{r}, b_{r0}, b_{r1})} \\
&= \dots
\end{align*}
To finish, observe that
\[
\D_{\geq r} / (a_{r}, b_{r0}, b_{r1}, \dots, b_{r,r-1}) =
\D_{\geq r+1}.
\]

Part (b) is similar; one just has to keep track of whether $\tau_{r}$
is present in $\dual{C}$, along with which powers of $\xi_{r}$ are
present.
\end{proof}

\begin{lemma}\label{lemma-bousfield-big}
Fix $q>R \geq 1$. 
\begin{enumerate}
\item In $\stable{A}$ we have
\begin{align*}
\bousfield{\H{D}} &= \bousfield{\H{D}_{\geq R+1}} \vee
\bigvee_{r=0}^{R} \left( \bousfield{a_{r}^{-1} \H{D}_{\geq r}} \vee
\bigvee_{i=0}^{r-1} \bousfield{b_{ri}^{-1}
\H{D}_{\geq r}/(a_{r}, b_{r0}, \dots, b_{r,i-1})} \right) \\
&= \bousfield{\H{D}_{\geq R+1}} \vee
\bigvee_{r=0}^{R}\left( \bousfield{a_{r}^{-1} \H{D}_{\geq r, 0, q}} \vee
\bigvee_{i=0}^{r-1} \bousfield{b_{ri}^{-1}
\H{D}_{\geq r, i+1, q}/(a_{r}, b_{r0}, \dots, b_{r,i-1})} \right).
\end{align*}
\item Given $\alg{C} \subseteq \A$ with profile function $((c_{i});
(e_{j}))$, in $\stable{C}$ we have
\begin{align*}
\bousfield{H(\alg{C} \cap \D)} &= \bousfield{H(\alg{C} \cap \D_{\geq
R+1})} \\
& \quad \vee
\bigvee_{r=0}^{R} \left( \bousfield{\tilde{a}_{r}^{-1} H(\alg{C} \cap \D_{\geq r})} \vee
\bigvee_{i=0}^{\min (c_{r}, r)-1} \bousfield{b_{ri}^{-1}
H(\alg{C} \cap \D_{\geq r})/(\tilde{a}_{r}, b_{r0}, \dots, b_{r,i-1})} \right) \\
&= \bousfield{H(\alg{C} \cap \D_{\geq R+1})} \\
& \quad \vee
\bigvee_{r=0}^{R} \left( \bousfield{\tilde{a}_{r}^{-1} H(\alg{C} \cap \D_{\geq r, 0, q})} \vee
\bigvee_{i=0}^{\min (c_{r}, r)-1} \bousfield{b_{ri}^{-1}
H(\alg{C} \cap \D_{\geq r, i+1, q})/(\tilde{a}_{r}, b_{r0}, \dots, b_{r,i-1})} \right),
\end{align*}
where $\tilde{a}_{r} = a_{r}$ if $\tau_{r} \neq 0$ in $\dual{C}$,
$\tilde{a}_{r}=0$ if $\tau_{r}=0$ in $\dual{C}$. (That is, when
$\tau_{r}=0$ in $\dual{C}$, omit terms involving $\tilde{a}_{r}^{-1}$,
and omit mention of $\tilde{a}_{r}$ in other terms.)
\end{enumerate}
\end{lemma}

This is an odd primary version of \cite[Lemma
5.3.3]{palmieri;steenrod-stable}.

\begin{proof}
We will prove part (a). Part (b) is proved almost identically; the
only difference is that if some $\tau_{r}$ or $\xi_{r}^{p^{i}}$ is
zero in $\dual{C}$, then we omit the corresponding self-maps and wedge
summands.

The first equality follows from \autoref{lemma-bousfield-hdr}. For
the second equality, we appeal to \autoref{cor-bousfield} and
the Nilpotence \autoref{lemma-nilpotence}: if $k\geq r+i+1$ and if
$\xi_{i}^{p^{k}}=0$ while $\xi_{n}^{p^{k}} \neq 0$, then
$b_{ri}b_{nk}$ is nilpotent, so in the above equality of Bousfield
classes, as long as $k \geq r+i+1$, then $\bousfield{b_{ri}^{-1}
\H{B}} = \bousfield{b_{ri}^{-1} \H{B}/b_{nk}}$. We apply this
repeatedly to the ``boxes'' $\xi_{n}^{p^{k}}$ which form the
difference between $\D_{\geq r}$ and $\D_{\geq r, i+1, q}$, to get
\begin{align*}
\bousfield{b_{ri}^{-1} H\D_{\geq r}} &=
\bousfield{b_{ri}^{-1} H\D_{\geq r}/b_{r+i+2,r+i+1}} \\
&=
\bousfield{b_{ri}^{-1} H\D_{\geq r}/(b_{r+i+2,r+i+1}, b_{r+i+3,r+i+2},
b_{r+i+3,r+i+1})} \\
&= \dots = \bousfield{b_{ri}^{-1} H\D_{\geq r, i+1, q}}
\end{align*}
for any $q$.  To allow a uniform choice of $q$ independent of $r$, we
choose $q > R$.

Similarly, since $a_{r} b_{nk}$ is nilpotent if $k \geq r$, then we
see that $\bousfield{a_{r}^{-1} \H{D}_{\geq r}} = \bousfield{a_{r}^{-1}
H\D_{\geq r,0,q}}$ for any $q>R$.
\end{proof}

\begin{proof}[Proof of \autoref{thm-bdd-nilpotence}] 
For part (a), as in the $p=2$ case (the proof of \cite[Theorem
5.1.6]{palmieri;steenrod-stable}) we can reduce to the smash
nilpotence case, and then we can reduce to a Bousfield class
argument. Summarizing the first paragraph of the $p=2$ proof, we end
up with an object $X^{\infty}$, and we want to know whether $H(\alg{C}
\cap \D) \Smash X^{\infty} = 0$, and we may assume that $\H{B}
\Smash X^{\infty} = 0$ for every sub-Hopf algebra $\alg{B} \subseteq
\alg{C} \cap \D$ with bounded profile function.

Since $\Fp$ is a (trivial) Hopf algebra with a bounded profile
function, we know that $H\Fp \Smash X^{\infty} = 0$. Thus any
map
\[
S^{0} \to \varinjlim_{R} H(\alg{C} \cap \D_{\geq R+1}) \Smash
X^{\infty} = H\Fp \Smash X^{\infty}
\]
is null; we have applied \autoref{lemma-colimit-hdr} here. Since
homotopy commutes with colimits, then for some $R$, a particular map
$S^{0} \to H(\alg{C} \cap \D_{\geq R+1}) \Smash X^{\infty}$ is null,
and by arguing as in the $p=2$ proof, this means that $H(\alg{C} \cap
\D_{\geq R+1}) \Smash X^{\infty} = 0$. Similarly, since each Hopf
algebra $\alg{C} \cap \alg{B}_{r,i+1}$ has bounded profile function,
we know that $H(\alg{C} \cap \alg{B}_{r,i+1}) \Smash X^{\infty} = 0$,
so we can repeat the colimit argument to see that
\begin{gather*}
b_{ri}^{-1} H(\alg{C} \cap \D_{\geq r, i+1, q})/(a_{r}, b_{r0}, \dots,
b_{r,i-1})  \Smash X^{\infty} = 0, \\
a_{r}^{-1} H\D_{\geq r, 0, q} \Smash X^{\infty} = 0,
\end{gather*}+
for all $r, i$. So by \autoref{lemma-bousfield-big}, we see that
$H(\alg{C} \cap \D) \Smash X^{\infty} = 0$, as desired.
\end{proof}

Presumably a classification of quasi-elementary sub-Hopf algebras of
$\A$ as in \autoref{conj-quasi} (see also
\autoref{remark-quasi-A}) would allow one to establish Bousfield
equivalences between $\alg{B}_{r}$ and some quasi-elementary $Q$, and
so we hope that we could then prove that the quasi-elementary sub-Hopf
algebras detect nilpotence.


\section{Discussion}\label{sec-discussion}

\subsection{Quasi-elementary sub-Hopf algebras of \texorpdfstring{$\P$}{P}}

One of main open questions here is the classification of
quasi-elementary sub-Hopf algebras of $\A$. Let's focus on
quasi-elementary sub-Hopf algebras of $\P$: if one could classify
those, it seems likely that this would carry over to a classification
of the quasi-elementary sub-Hopf algebras of $\A$: see
\autoref{remark-quasi-A}. Here is at least some of what we know
about the situation with $\P$:
\begin{enumerate}
\item Every elementary Hopf algebra is quasi-elementary. The maximal
elementary sub-Hopf algebras of $\A$ were described in
\autoref{subsec-hopf}, and it is easy to figure out the
elementary sub-Hopf algebras of $\P$ from there: the maximal ones have
profile function $\mathbf{e}(r) := (\underbrace{0, \dots, 0}_{r}, r+1,
r+1, \dots)$, and every list $(k_{1}, k_{2}, k_{3}, \dots)$ of
non-negative integers which is term-wise no larger than some
$\mathbf{e}(r)$ is a valid profile function and corresponds to an
elementary Hopf algebra.  That is, the profile function fits inside a
rectangle like this:

\begin{tikzpicture}[scale=0.7]
\draw (0,4.5) -- (0,0) -- (10,0);
\draw (4.5, 0.5) node {$\xi_{r}$};
\draw (4.5, 1.5) node {$\vdots$};
\draw (4.7, 2.5) node {$\xi_{r}^{p^{r-1}}$};
\draw[gray] (4.5, 3.5) node {$\xi_{r}^{p^{r}}$};
\draw (4, 0) -- (4, 3) -- (8, 3);
\end{tikzpicture}

\item Fix $r \geq 2$. The computation given in
\cite[6.3]{wilkerson}, and see also \cite[Theorem 1.1]{lin;jpaa-I},
generalizes to show that Hopf algebras with profile function
\[
( \underbrace{0, \dots, 0}_{r-2}, 1, n_{r}, n_{r+1}, \dots )
\]
with each $n_{i} \leq r$ are quasi-elementary: the profile function
fits inside

\begin{tikzpicture}[scale=0.7]
\draw (0,4.5) -- (0,0) -- (10,0);
\draw (3.6, 0.5) node {$\xi_{r-1}$};
\draw (4.5, 0.5) node {$\xi_{r}$};
\draw (4.5, 1.5) node {$\vdots$};
\draw (4.7, 2.5) node {$\xi_{r}^{p^{r-1}}$};
\draw[gray] (4.5, 3.5) node {$\xi_{r}^{p^{r}}$};
\draw (3, 0) -- (3,1) -- (4,1) -- (4, 3) -- (8, 3);
\end{tikzpicture}

These Hopf algebras include the ones in (a).
\end{enumerate}

\begin{conjecture}\label{conj-quasi}
The quasi-elementary sub-Hopf algebras of $\P$ are those described in
(b).
\end{conjecture}

(This was stated as a theorem in \cite{nakano-palmieri}.)

Note that one can compute the cohomology of a Hopf algebra as in (b)
--- call it $\alg{Q}$ --- by imitating Wilkerson's work in
\cite[6.3]{wilkerson}. There is a spectral sequence converging to the
cohomology with $E_{2}$-page of the form
\[
\exterior (h_{ts} : \xi_{t}^{p^{s}} \neq 0 \text{ in $\dual{Q}$}) \otimes \Fp
[b_{ts} : \xi_{t}^{p^{s}} \neq 0 \text{ in $\dual{Q}$}].
\]
The differentials only affect the exterior part, so the cohomology
ends up being the tensor product of the polynomial part with
a subquotient of the exterior part. That is, there is an
$F$-isomorphism
\[
H^{*}(\alg{Q}) \xrightarrow{F-\text{iso}} \Fp [b_{ts} :
\xi_{t}^{p^{s}} \neq 0 \text{ in $\dual{Q}$}].
\]
Any inclusion $\alg{Q} \hookrightarrow \alg{Q}'$ induces the obvious
map on the polynomial part of the cohomology, so one can then compute
the limit of the cohomology algebras of these, at least up to
$F$-isomorphism. So if the conjecture holds, one would get a close
analog to what happens at the prime 2, computationally as well as
theoretically. In particular, one should be able to deduce the
following from the previous conjecture:

\begin{conjecture}\label{conj-cohomology}
Let $S = \Fp [b_{ts}: s<t]$. Let $J \subseteq S$
be the ideal
\[
J = 
(b_{ts} b_{vu} \,:\, u \geq t \text{ unless $s=0$ and $u=t$}).
\]
Then there is an $F$-isomorphism $\lim H^{*}(\alg{Q}) \to S/J$,
where the limit is taken over the quasi-elementary sub-Hopf algebras
of $\P$. There is an action of $\P$ on $S/J$; this $F$-isomorphism
is a maps of $\P$-algebras, and the induced map on the invariant
subalgebras is again an $F$-isomorphism.
\end{conjecture}




Combined with \autoref{conj-f-iso-quasi}, one should be able to see
that there is an $F$-isomorphism
\[
H^{*}(\D \cap \P) \to \lim H^{*}(\alg{Q}),
\]
which is a map of $\P$-algebras. Compare this to the $p=2$
computations of $H^{*}(\D)$ and $\lim H^{*}(\alg{E})$ up to
$F$-isomorphism in \cite[Proposition 6.5.1]{palmieri;steenrod-stable}.

Here is some of what we know about the status of
\autoref{conj-quasi}.  First, if for some $r$ the profile
function of a Hopf algebra $\alg{B}$ has the form
\[
( \underbrace{0, \dots, 0}_{r-2}, 1, b_{r}, b_{r+1}, \dots )
\]
and $b_{i} > r$ for some $i$, then $\alg{B}$ cannot be
quasi-elementary by \cref{lemma-nilpotence,lemma-quasi}. As a result,
if $1$ is the first nonzero entry in the profile function, then the
Hopf algebra is quasi-elementary if and only if it is of the form
given in (b).  So to address the classification question, we need to
consider Hopf algebras where the first nonzero entry is larger than
$1$. It turns out --- see the paragraph after
\autoref{question-lin} --- that the following case seems to be
difficult: the case when the profile function looks like
\[
(0, \dots, 0, \underset{t}{s}, 0, \dots, 0, \underset{n}{s+t}, \dots),
\]
with $1<s<t$, so for example $(0,2,0,4,0,2)$.  The author is not able
to determine whether the corresponding Hopf algebra is
quasi-elementary. If it is not, it lends support to the claimed
classification. If it is, then note that it has a sub-Hopf algebra
with profile function $(0,1,0,4,0,2)$ which is not quasi-elementary,
and it would perhaps be surprising to have $\alg{B} \subseteq \alg{Q}
\subseteq \alg{P}$ with $\alg{Q}$ quasi-elementary and $\alg{B}$ not
quasi-elementary.


\subsection{The Nilpotence Lemma}

As already discussed, the Nilpotence \autoref{lemma-nilpotence} is an
odd primary analog of a result of Lin's \cite[Theorem
3.1]{lin;jpaa-I}.  If one frames Lin's result the right way, it can be
stated as an if-and-only-if statement. Can one do the same thing here?

To set this up, let's fix integers $s,t,n,k$ with $s<t$ and $k<n$ ---
we want to consider whether $b_{ts} b_{nk}$ is nilpotent in the
cohomology of an appropriate Hopf algebra, and by
\autoref{thm-nilpotence} we know that $b_{ts}$ is nilpotent if $s \geq
t$, and similarly for $b_{nk}$. Let's also assume that $t<n$. To keep
things simple, let's work with a Hopf algebra $\alg{C}$ whose profile
function looks like
\[
(0, \dots,0, \underset{t}{s}, 0, \dots, 0, \underset{n}{c_{n}}, \dots)
\]
where $c_{n} > k$. Therefore $\xi_{t}^{p^{s}}$ and $\xi_{n}^{p^{k}}$
are nonzero and primitive in $\dual{C}$.

At the prime 2, Lin's theorem says that in this situation, $h_{nk}
h_{ts}$ is nilpotent if $k \geq t$, while if $k<t$, by restricting to
an appropriate elementary sub-Hopf algebra one can see that $h_{nk}
h_{ts}$ is non-nilpotent. That is, with $s$, $t$, and $n$ fixed,
$h_{ts} h_{nk}$ is nilpotent if and only if $k \geq t$.

When $p$ is odd, here are two cases to consider:
\begin{enumerate}
\item $c_{n}= \infty$: no power of $\xi_{n}$ is zero in $\dual{C}$,
\item $c_{n}>k$: $\xi_{n}^{p^{k}} \neq 0$ while some power of
$\xi_{n}$ is zero.
\end{enumerate}
Unfortunately, in case (a) the author can say nothing beyond what is
already in the lemma. For case (b), though, we can say more.  We have
the reduced coproduct $\xi_{n+t} \mapsto \xi_{n}^{p^{t}} \otimes
\xi_{t}$, so we get the relation $h_{nt} h_{t0} = 0$ in cohomology. If
we apply Steenrod operations, we get more relations:
\begin{align*}
\beta \tildePP{0}: & \quad b_{nt} h_{t1} = h_{n,t+1} b_{t0}, \\
\beta \tildePP{1}: & \quad b_{nt}^{p} b_{t1} = b_{n,t+1} b_{t0}^{p},\\
\tildePP{p}: & \quad b_{nt}^{p^{2}} b_{t2} = b_{n,t+2} b_{t0}^{p^{2}},\\
\intertext{and then apply $\tildePP{p^{2}}$, $\tildePP{p^{3}}$, etc. so
in general for $d \geq 1$ we get}
& \quad b_{nt}^{p^{d}} b_{td} = b_{n,t+d} b_{t0}^{p^{d}}.
\intertext{Then with $i \geq 0$ we apply $(\tildePP{0})^{i}$ to get}
& \quad b_{n,t+i}^{p^{d}} b_{t,d+i} = b_{n,t+d+i} b_{t,i}^{p^{d}}.
\intertext{If we let $s=d+i$, this becomes}
& \quad b_{n,s+t-d}^{p^{d}} b_{ts} = b_{n,s+t} b_{t,s-d}^{p^{d}}.
\end{align*}
If $\xi_{n}^{p^{s+t}}=0$, then the right side is zero, so we see that
$b_{n,s+t-d} b_{ts}$ is nilpotent whenever $s \geq d \geq 1$. Letting
$k=s+t-d$, we conclude that if $\xi_{n}^{p^{s+t}}=0$, then
$b_{nk}b_{ts}$ is nilpotent whenever $t \leq k \leq s+t-1$. The
Nilpotence Lemma tells us that if $k \geq s+t+1$ (and $\xi_{n}^{p^{k}}
\neq 0$ is nonzero), then $b_{nk} b_{ts}$ is nilpotent. This leaves
one open case: what happens if $k=s+t$? Wilkerson's result
\cite[6.3]{wilkerson} tells us what happens when $s=0$: $b_{nt}
b_{t0}$ need not be nilpotent. What happens if $s>0$?

\begin{question}\label{question-lin}
Suppose that the sub-Hopf algebra $\alg{C}$ of $\P$ has a finite
profile function, and assume that $\xi_{t}^{p^{s}}$ and
$\xi_{n}^{p^{s+t}}$ are nonzero and primitive in $\dual{C}$. Suppose
that $0<s<t$ and $s+t<n$. Suppose further that $\xi_{t}^{p^{s+1}}=0$ and
$\xi_{n}^{p^{s+t+1}}=0$. Is $b_{ts} b_{n,s+t}$ nilpotent?
\end{question}

For example, is $b_{21} b_{43}$ nilpotent in the Hopf algebra with
profile function $(0,2,0,4,0,2)$?  Answering this should determine
whether this Hopf algebra is quasi-elementary.  Completely answering
the question would help to classify the quasi-elementary sub-Hopf
algebras of $\P$.

Note that the situation with products of the form $a_{t} b_{nk}$ is
simpler. Assume that $k<n$. Then $a_{t} b_{nk}$ is nilpotent if and
only if $k \geq t$: the $k \geq t$ case comes from the Nilpotence
Lemma, while if $k<t$, then $a_{t} b_{nk}$ is detected in the
cohomology of an elementary sub-Hopf algebra of $\A$ and hence is
non-nilpotent.

\begin{remark}\label{remark-quasi-A}
Therefore as far as quasi-elementary sub-Hopf algebras of $\A$ go, if
we can understand the quasi-elementary sub-Hopf algebras of $\P$, we
should also be able to understand those in $\A$. If
\autoref{conj-quasi} is correct, for example, then the maximal
quasi-elementary sub-Hopf algebras of $\A$ should be the ones with
profile function
\[
(( \underbrace{0, \dots, 0}_{r-1}, 1, r+1, r+1, r+1, \dots ); 
(\underbrace{1, \dots, 1}_{r+1}, 2, 2, 2, \dots))
\]
for some $r \geq 1$, along with the elementary Hopf algebras
$\alg{E}(-1)$ and $\alg{E}(0)$ from \autoref{subsec-hopf}.
\end{remark}


\subsection{Computations in \texorpdfstring{$H^{*}(\D)$}{H*(D)}}

We note that $\D$ is the sub-Hopf algebra generated by the elementary
sub-Hopf algebras of $\A$, and also generated by the quasi-elementary
sub-Hopf algebras of $\A$. It is also the sub-Hopf algebra generated
by the Margolis elements $\pst{s}{t}$ with $s<t$ along with the Milnor
primitives $Q_{n}$, and those elements are precisely the ones which
determine projectivity and vanishing lines in cohomology
\cite{miller-wilkerson}. So it is natural that $\D$ is so prominent in
this work.

It would be good to understand something about the cohomology of
$H^{*}(\D)$ in order to provide examples and applications of
\cref{thm-nilpotence,thm-f-iso}. We recall an old
observation, that for a graded connected Hopf algebra $\alg{B}$, the
vector space $H^{1}(\alg{B})$ is isomorphic to the primitives in
$\dual{B}$, with a primitive $x \in \dual{B}$ corresponding to the
cohomology class represented by $[x]$ in the cobar complex. This was
originally probably a folk theorem. It can be deduced by dualizing the
results from Wall \cite[p.\@ 442]{wall;steenrod}, together with Milnor
and Moore's observation \cite[Proposition 3.10]{milnor-moore} that the
indecomposables in $\alg{B}$ are dual to the primitives in $\dual{B}$.

\begin{proposition}\label{prop-D-primitives}
$\xi_{t}^{p^{s}}$ is primitive in $\dual{D}$ if and only if $s
\geq \floor{t/3}$. $\tau_{n}$ is primitive in $\dual{D}$ if and only
if $n = 0$.
\end{proposition}

The condition on $\xi_{t}^{p^{s}}$ being primitive is the same as
given by Harker in the $p=2$ case \cite[2.7--2.8]{harker}.

\begin{proof}
This is just a computation in $\dual{D} = \A / (\xi_{k}^{p^{k}}: k
\geq 1)$. The primitivity of $\tau_{n}$ if and only if $n=0$ is clear.
It remains to show that $\xi_{3n}^{p^{s}}$, $\xi_{3n+1}^{p^{s}}$, and
$\xi_{3n+2}^{p^{s}}$ are primitive in $\dual{D}$ if and only if $s
\geq n$. There are three cases and they are very similar, so we do
the last one. The reduced coproduct on $\xi_{3n+2}$ is
\[
\sum_{j=1}^{3n+1} \xi_{3n+2-j}^{p^{j}} \otimes \xi_{j}.
\]
We need to show that the coproduct on $\xi_{3n+2}^{p^{n-1}}$ is
nonzero, while the coproduct on $\xi_{3n+2}^{p^{n}}$ is zero. If $j=n,
n+1$, we get the terms
\[
\xi_{2n+2}^{p^{n}} \otimes \xi_{n} + \xi_{2n+1}^{p^{n+1}} \otimes \xi_{n+1}.
\]
If we raise these to the power $p^{n-1}$, the terms are nonzero:
\[
\xi_{2n+2}^{p^{2n-1}} \otimes \xi_{n}^{p^{n-1}} + \xi_{2n+1}^{p^{2n}} \otimes \xi_{n+1}^{p^{n-1}}.
\]
Therefore $\xi_{3n+2}^{p^{n-1}}$ is not primitive.

Now consider the coproduct on $\xi_{3n+2}^{p^{n}}$: the $j$th term is
\[
\xi_{3n+2-j}^{p^{j+n}} \otimes \xi_{j}^{p^{n}}.
\]
If $j \leq n$, then the second factor is zero. If $j \geq n+1$, then
the first factor is zero: if $j \geq n+1$, then $3n+2-j\leq 2n+1$
while $j+n \geq 2n+1$.

The computations for $\xi_{3n}$ and $\xi_{3n+1}$ are similar and
a bit simpler than this one.
\end{proof}

\begin{corollary}
$H^{1}(\D)$ is the vector space spanned by $\{ h_{ts} \,:\, s \geq
\floor{t/3} \} \cup \{ a_{0} \}$.
\end{corollary}

The author believes that the following is true, and we provide a
sketch of a proof. Some details are lacking, though.

\begin{guess}
For each $n \geq 1$, the element $b_{n,n-1} \in H^{2}(\D)$ is
invariant under the $\A$-action. Therefore some $p$th
power of $b_{n,n-1}$ lifts to $H^{*}(\A)$.
\end{guess}

\begin{proof}[Sketch of proof]
We first note that $b_{n,n-1} \in H^{2}(\D)$ is detected by the
restriction map to $H^{*}(D(\pst{n-1}{n}))$, and so is non-nilpotent.

The action of $\A$ on $H^{*}(\D)$ comes from the action of $\A$ on
$\D$ itself. It is convenient to recast this as a coaction of
$\dual{A}$ on $H^{*}(\D)$.  As at the prime 2 --- see Lemmas A.3 and
A.5 of \cite{palmieri;nilp2} --- the coaction on an element $h_{ts}$
comes from the coproduct on $\xi_{t}^{p^{s}}$ and so has the form $1
\otimes h_{ts} + \sum_{i,j} x \otimes h_{ij}$ where $x \in \A$ and $i
< t$, $j \geq s$. In particular, the elements $h_{n,n-1}$ are
invariant.  We would like at this point to appeal to the naturality of
May's Steenrod operations, but it is not clear how to do
that. Naturality should imply that this coaction commutes with the
Steenrod operations on Hopf algebra cohomology, so the same sort of
coaction formula should hold for $b_{ts}$, and therefore $b_{n,n-1}$
should be invariant for each $n$. That is, the subalgebra of
$\A$-invariants in $H^{*}(\D)$ should be closed under the action of
May's Steenrod operations, so since each $h_{n,n-1}$ is invariant,
each $b_{n,n-1}$ should also be.
\end{proof}

In addition, if \autoref{conj-cohomology} is correct (see also
\autoref{remark-quasi-A}), then for any $n \geq 2$ one would get
invariant monomials of the form
\[
b_{n-1,0}^{i} b_{n,1}^{j_{1}} b_{n,2}^{j_{2}} \dots b_{n,n-1}^{j_{n-1}},
\]
as long as $j_{n-1}>0$. Indeed, rather than trying to compute the
cohomology of $\D$, the path forward should probably be to classify
the quasi-elementary sub-Hopf algebras and then to work with their
cohomology.


\subsection{Other results}

If one would like results that do not rely on the unpublished work
\cite{burklund-levy-pstragowski}, then we note that this was only used
in \autoref{sec-binom}. In particular nothing in
\autoref{sec-D} relies on it, so we have:

\begin{theorem}
Suppose that $\alg{C} \subseteq \A$ has a finite profile
function. Then
\begin{enumerate}
\item The sub-Hopf algebra $\alg{C} \cap \D$ detects nilpotence in
$\stable{C}$.
\item The map
\[
\overline{\imath}^{*}: \Ext_{\alg{C}}^{*}(\Fp, \Fp) \to \left( \Ext_{\alg{C}
\cap \D}^{*}(\Fp, \Fp)\right)^{\alg{C} \modmod \alg{C} \cap \D}
\]
is an $F$-isomorphism.
\item The sub-Hopf algebras of $\alg{C}\cap \D$ with bounded profile
function detect nilpotence in $\stable{C}$.
\end{enumerate}
\end{theorem}

We prove the Burklund-Levy-Pstr\k{a}gowski theorem when $p=3$ in the
appendix, so we have:

\begin{theorem}
Fix $p=3$ and suppose that $\alg{C}$ is a sub-Hopf algebra of $\A$.
\begin{enumerate}
\item The sub-Hopf algebra $\alg{C} \cap \D$ detects nilpotence in
$\stable{C}$.
\item The map
\[
\overline{\imath}^{*}: \Ext_{\alg{C}}^{*}(\Fp, \Fp) \to \left( \Ext_{\alg{C}
\cap \D}^{*}(\Fp, \Fp)\right)^{\alg{C} \modmod \alg{C} \cap \D}
\]
is an $F$-isomorphism. 
\item The sub-Hopf algebras of $\alg{C} \cap \D$ with bounded profile
function detect nilpotence in $\stable{C}$.
\end{enumerate}
\end{theorem}


\appendix
\section{Nilpotence of \texorpdfstring{$b_{nn}$}{bnn} when 
\texorpdfstring{$p=3$}{p=3}}\label{sec-bnn-nilp}

We fix $p=3$ for this section. The goal is to provide a proof of the
nilpotence of $b_{nn}$ in $H^{*}(\P_{\geq n})$, independent of
\autoref{thm-blp}.

We start by mentioning that Nakamura proved that $b_{11} \in
H^{*}(\P)$ is nilpotent. More precisely, he showed \cite[Proposition
1.1c)]{nakamura} that $b_{11}^{2} = h_{11} x$ for some $x$, and since
$h_{11}^{2}=0$, one immediately deduces that $b_{11}^{4}=0$. His proof
is a calculation in the cobar complex, and it is not easily
generalizable (at least not by the author) to either higher primes or
to $b_{nn}$ for larger values of $n$. His element $x$ looks like a
Massey product, though, and that observation inspired the proof of the
following result.

\begin{theorem}\label{thm-3}
If $p=3$, then $b_{nn} \in H^{2}(\P_{\geq n})$ is nilpotent.
\end{theorem}

The proof here follows the sketch in \cite[Section
B.3.1]{palmieri;steenrod-stable}.

\begin{proof}
We work with $\P_{\geq n}$ and its cohomology.  The reduced coproduct
$\xi_{2n} \mapsto \xi_{n}^{p^{n}} \otimes \xi_{n}$ yields the relation
\[
h_{nn} h_{n0} = 0.
\]
If we apply Steenrod operations to this, we get the following, where
each line is labeled by the Steenrod operation which has been applied:
\begin{align*}
\beta \tildePP{0}: b_{nn} h_{n1} - h_{n,n+1} b_{n0} &= 0, \\
\tildePP{1}: b_{nn}^{p} h_{n2} - h_{n,n+2} b_{n0}^{p} &= 0, \\
\tildePP{p}: b_{nn}^{p^{2}} h_{n3} - h_{n,n+3} b_{n0}^{p^{2}} &= 0 \\
\vdots \\
\tildePP{p^{k-2}}: b_{nn}^{p^{k-1}} h_{nk} - h_{n,n+k} b_{n0}^{p^{k-1}} &= 0
\end{align*}
That is, $b_{nn}^{p^{k-1}} h_{nk} = h_{n,n+k} b_{n0}^{p^{k-1}}$ for
any $k \geq 1$. So letting $k=n$, we have
\begin{equation}\label{eqn-massey-setup1}
h_{nn} b_{nn}^{p^{n-1}} = h_{n,2n} b_{10}^{p^{n-1}},
\end{equation}
and letting $k=2n$, we have
\[
b_{nn}^{p^{2n-1}} h_{n,2n} = h_{n,3n} b_{n0}^{p^{2n-1}}.
\]
We apply both of these relations to get
\begin{align}\label{eqn-massey-setup2}
h_{nn} b_{nn}^{p^{n-1} + p^{2n-1}} &= 
h_{n,2n} b_{10}^{p^{n-1}} b_{nn}^{p^{2n-1}} \nonumber \\
&= h_{n,3n} b_{10}^{p^{n-1} + p^{2n-1}}.
\end{align}
(More generally, $h_{nn}b_{nn}^{M} = h_{n,mn} b_{n0}^{M}$, where
$M=p^{n-1} + p^{2n-1} + \dots + p^{(m-1)n-1}$, but we don't need
this.)

Now we have the Massey product $\massey{h_{n,3n}, h_{n,2n}, h_{nn}}$,
defined because $h_{n,n+k} h_{nk} = 0$ for all $k$: apply the algebra
map $(\tildePP{0})^{k}$ to the relation $h_{nn} h_{n0} = 0$. We have the
reduced coproduct 
\begin{align}\label{eqn-A}
\xi_{3n} \mapsto & \sum_{i=1}^{3n-1} \xi_{3n-i}^{p^{i}} \xi_{i} \nonumber \\
&= \sum_{i=n}^{2n} \xi_{3n-i}^{p^{i}} \xi_{i} \nonumber \\
&= \xi_{2n}^{p^{n}} \xi_{n} + \xi_{n}^{p^{2n}} \otimes \xi_{2n} + 
\sum_{i=n+1}^{2n-1} \xi_{3n-i}^{p^{i}} \xi_{i}.
\end{align}
We note that if $n \leq i \leq 2n-1$, then $\xi_{i}$ is primitive in
$\dual{P}_{\geq n}$. Hence the tensor factors in the sum give
cohomology classes, and the above coproduct on $\xi_{3n}$ gives us
this relation in cohomology:
\[
0 = \massey{h_{n,2n}, h_{nn}, h_{n0}} + \sum_{i=n+1}^{2n-1}
h_{3n-i,i} h_{i,0}.
\]
Apply $\tildePP{0}$ to get
\[
\massey{h_{n,3n}, h_{n,2n}, h_{nn}} = - \sum_{i=n+1}^{2n-1}
h_{3n-i,i+1} h_{i,1}.
\]
Note that the right side is a product of exterior elements, and so it
is nilpotent; thus $\massey{h_{n,3n}, h_{n,2n}, h_{nn}}$ is nilpotent.

Multiply this Massey product by $b_{n0}^{N}$ where $N=p^{2n-1} +
p^{n-1} + p^{n-1}$: first we get
\[
b_{n0}^{N} \massey{h_{n,3n}, h_{n,2n}, h_{nn}}
\in \massey{h_{n,3n} b_{n0}^{p^{2n-1} + p^{n-1}}, h_{n,2n}
b_{n0}^{p^{n-1}}, h_{nn}}.
\]
The Massey product on the left side consists of a single element, for
degree reasons. The Massey product on the right is a set, containing
the single element on the left. Now we apply \eqref{eqn-massey-setup1}
and \eqref{eqn-massey-setup2}:
\begin{align*}
b_{n0}^{N} \massey{h_{n,3n}, h_{n,2n}, h_{nn}}
&\in \massey{h_{n,3n} b_{n0}^{p^{2n-1} + p^{n-1}}, h_{n,2n}
b_{n0}^{p^{n-1}}, h_{nn}} \\
&= \massey{h_{nn} b_{nn}^{p^{2n-1} + p^{n-1}}, h_{nn}
b_{nn}^{p^{n-1}}, h_{nn}} \\
&\ni b_{nn}^{N} \massey{h_{nn}, h_{nn}, h_{nn}} = b_{nn}^{N+1}.
\end{align*}
Therefore 
\[
b_{nn}^{N+1} - b_{n0}^{N}\massey{h_{n,3n}, h_{n,2n} h_{nn}}
\]
is in the indeterminacy of the Massey product
\[
\massey{h_{nn} b_{nn}^{p^{2n-1} + p^{n-1}}, h_{nn}
b_{nn}^{p^{n-1}}, h_{nn}}.
\]
We claim that every element in the indeterminacy of this Massey
product is nilpotent.

\begin{lemma}\label{lemma-massey-indeterminacy} 
Suppose that $\Omega$ is a differential graded algebra over a field
whose characteristic is not two, and assume that the cohomology
algebra $H^{*}\Omega$ is graded commutative.  If $x$ and $z$ are
odd-dimensional classes in $H^{*}\Omega$ and if $\massey{x,y,z}$ is
defined, then every element in the indeterminacy of $\massey{x,y,z}$
is nilpotent.
\end{lemma}

\begin{proof}
The indeterminacy of a length 3 Massey product has an explicit
formula, given by May in \cite{may;massey}: every element in the
indeterminacy of $\massey{x,y,z}$ has the form $x a + b z$, where $a$
and $b$ are cohomology classes in the appropriate dimensions. Since $x$
and $z$ are odd-dimensional cohomology classes in a graded commutative
algebra and since we are not working in characteristic two, both $x$
and $z$ are nilpotent and hence so is every element in the
indeterminacy.
\end{proof}

Since $\massey{h_{n,3n}, h_{n,2n} h_{nn}}$ is nilpotent, we see that
$b_{nn}^{N+1}$ is the sum of two nilpotent elements, and hence is
nilpotent itself.
\end{proof}

This proof does not easily generalize to larger primes for two
reasons: first, at larger primes the various terms in \eqref{eqn-A}
are not all primitive and so do not correspond to cohomology
classes. Second, \autoref{lemma-massey-indeterminacy} is not known
to hold at larger primes. See \cite[Conjectures
B.3.9--10]{palmieri;steenrod-stable} for a possible replacement for
the lemma.


\def\cftil#1{\ifmmode\setbox7\hbox{$\accent"5E#1$}\else
  \setbox7\hbox{\accent"5E#1}\penalty 10000\relax\fi\raise 1\ht7
  \hbox{\lower1.15ex\hbox to 1\wd7{\hss\accent"7E\hss}}\penalty 10000
  \hskip-1\wd7\penalty 10000\box7} \def\cprime{$'$}
\providecommand{\bysame}{\leavevmode\hbox to3em{\hrulefill}\thinspace}
\providecommand{\MR}{\relax\ifhmode\unskip\space\fi MR }
\providecommand{\MRhref}[2]{%
  \href{http://www.ams.org/mathscinet-getitem?mr=#1}{#2}
}
\providecommand{\href}[2]{#2}

\end{document}